\documentclass[reqno,11pt]{amsart}
\usepackage{graphicx}
\usepackage{amsfonts,amsmath, amssymb}
\usepackage{marginnote}
\usepackage{color}
\usepackage{mathrsfs}
\usepackage{tikz}
\numberwithin{equation}{section}

\newtheorem*{pf}{Proof}

\DeclareMathOperator\dif{d\!}

\newtheorem{theorem}{Theorem}[section]

\newtheorem{lemma}[theorem]{Lemma}

\newtheorem{remark}[theorem]{Remark}

\def\v{\varepsilon}

\def\f{\frac}

\begin{document}

 \title[BE with Large External Potential]{Global Stability of Boltzmann Equation with Large External Potential for a Class of Large Oscillation Data}

 \author{Guanfa Wang}
 \address{Institute of Applied Mathematics, Academy of Mathematics and Systems Science, Chinese Academy of Sciences, Beijing 100190, China, and University of Chinese Academy of Sciences}
 \email{wgf@amss.ac.cn}

 \author{Yong Wang}
 \address{Institute of Applied Mathematics, Academy of Mathematics and Systems Science, Chinese Academy of Sciences, Beijing 100190, China, and University of Chinese Academy of Sciences}
 \email{yongwang@amss.ac.cn}

 \begin{abstract}
 	
In this paper, we investigate the stability of Boltzmann equation with large external potential in $\mathbb{T}^3$. For a class of initial data with large oscillations in $L^\infty_{x,v}$ around the local Maxwellian, we prove the existence of a global solution to the Boltzmann equation provided the initial perturbation  is suitably small in $L^2$-norm. The large time behavior of the Boltzmann solution with exponential decay rate is  also obtained. This seems to be the first result on the perturbation theory of large-amplitude non-constant equilibriums for large-amplitude initial data. 

 \end{abstract}


 \date{\today}
 \maketitle


 \thispagestyle{empty}


 \section{Introduction}
  In this paper, we consider the Boltzmann equation with external potential
  \begin{equation}\label{eq:main}
   F_t +v \cdot \nabla_x F -\nabla\Phi(x) \cdot \nabla_v F =Q(F,F),
  \end{equation}
  supplemented with initial data
  \begin{equation}\label{initial data}
   F(0,x,v)=F_0(x,v).
  \end{equation}
  The unknown $F=F(t,x,v)\ge0$ represents for the density distribution function of gas particles with position $x\in\mathbb{T}^3$ and particle velocity $v\in\mathbb{R}^3$ at time $t>0$.
  The collision term $Q(F,F)$ is an integral with respect to velocity variable only, and it takes the non-symmetric bilinear form
  \begin{equation}\nonumber
   \begin{split}
    Q(F_1,F_2) &=\int_{\mathbb{R}^3} \int_{\mathbb{S}^2} B(v-u,\omega) F_1(u') F_2(v') \dif \omega \dif u -\int_{\mathbb{R}^3} \int_{\mathbb{S}^2} B(v-u,\omega) F_1(u) F_2(v) \dif \omega \dif u
    \\& \triangleq Q_+(F_1,F_2)-Q_-(F_1,F_2).
   \end{split}
  \end{equation}
  Here the post-collision velocity pair $(v',u')$ and the pre-collision velocity pair $(v,u)$ satisfy the relation
  \begin{equation}\nonumber
   u'+v'=u+v \quad\text{and}\quad |u'|^2+|v'|^2=|u|^2+|v|^2.
  \end{equation}
  The collision kernel $B(v-u,\omega)$ depends only on the relative velocity $|v-u|$ and $\cos\theta:=(v-u)\cdot\omega/|v-u|$, and it is assumed throughout the paper to take the form
  \begin{equation}\nonumber
   B(v-u,\omega)=|v-u|^\gamma b(\theta),
  \end{equation}
  with $b(\theta)$ satisfying the angular cutoff assumption
  \begin{equation}\nonumber
   0\le b(\theta)\le C |\cos\theta|.
  \end{equation}
  Throughout the paper, we consider the hard potential case, i.e. $0\le\gamma\le1$. The external potential $\Phi$ in \eqref{eq:main} is a given periodic function depending only on the spatial variable $x\in\mathbb{T}^3$. We assume the function $\Phi(x)\in C^3(\mathbb{T}^3)$ and denote
  \begin{equation}\label{def:the C3 norm of Phi}
   M\triangleq\|\Phi\|_{C^3(\mathbb{T}^3)}.
  \end{equation}
  Without loss of generality, we assume that $\Phi(x)\ge0$, otherwise one can replace $\Phi$ by $M+\Phi$ while the equation \eqref{eq:main} remains the same.

  For given external potential, it is direct to check that the local Maxwellian
  \begin{equation}\nonumber
   \mu_E(x,v)=\exp\left\{-\frac{|v|^2}2-\Phi(x)\right\}=\mu(v)e^{-\Phi(x)},
  \end{equation}
  is a steady solution to the Boltzmann equation \eqref{eq:main}. We define the perturbation function
  \begin{equation}\nonumber
   f(t,x,v):=\frac{F(t,x,v)-\mu_E(x,v)}{\sqrt{\mu_E(x,v)}},
  \end{equation}
  then the Boltzmann equation is rewritten as
  \begin{equation}\label{eq:equation of f}
   f_t+v\cdot\nabla_x f-\nabla\Phi\cdot\nabla_v f+e^{-\Phi}Lf =e^{-\frac\Phi2} \Gamma(f,f),
  \end{equation}
  where $L$ is the standard linearized operator given by
  \begin{equation}\nonumber
   Lf=-\frac1{\sqrt\mu} \left\{Q(\mu,\sqrt\mu f)+Q(\sqrt\mu f,\mu)\right\}=\nu(v)f-Kf,
  \end{equation}
  with $K=K_2-K_1$ defined by
  \begin{align}
   K_1 f(v)&=\int_{\mathbb{R}^3} \int_{\mathbb{S}^2} B(v-u,\omega) \sqrt{\mu(v)\mu(u)}f(u) \dif \omega \dif u,\label{def:K1}\\
   K_2 f(v)&=\int_{\mathbb{R}^3} \int_{\mathbb{S}^2} B(v-u,\omega) \sqrt{\mu(u)\mu(u')}f(v') \dif \omega \dif u\nonumber\\
   &\quad+\int_{\mathbb{R}^3} \int_{\mathbb{S}^2} B(v-u,\omega) \sqrt{\mu(u)\mu(v')}f(u') \dif \omega \dif u,\label{def:K2}
  \end{align}
  and the collision frequency
  \begin{equation}\nonumber
   \nu(v)=\int_{\mathbb{R}^3} \int_{\mathbb{S}^2} |v-u|^\gamma \mu(u) b(\theta) \dif \omega \dif u \sim (1+|v|)^\gamma,\quad 0\le\gamma\le1.
  \end{equation}
  The nonlinear term $\Gamma(f,f)$ is given by
  \begin{equation}\nonumber
   \begin{split}
    \Gamma(f,f)&:=\frac1{\sqrt\mu} Q_+(\sqrt\mu f,\sqrt\mu f)-\frac1{\sqrt\mu} Q_-(\sqrt\mu f,\sqrt\mu f)
    \\&=:\Gamma_+(f,f)-\Gamma_-(f,f).
   \end{split}
  \end{equation}

  Let $F$ be a solution of the Boltzmann equation \eqref{eq:main}. Since an external potential $\Phi$ is included, generally we only have the conservation of mass and energy:
  \begin{align}
   &\int_{\mathbb{T}^3} \int_{\mathbb{R}^3} \left\{F(t,x,v)-\mu_E(x,v)\right\} \dif v \dif x=\int_{\mathbb{T}^3} \int_{\mathbb{R}^3} \left\{F_0(x,v)-\mu_E(x,v)\right\} \dif v \dif x=M_0,\label{eq:conservation of mass}\\
   &\int_{\mathbb{T}^3} \int_{\mathbb{R}^3} \left[\frac{|v|^2}2+\Phi(x)\right] \left\{F(t,x,v)-\mu_E(x,v)\right\} \dif v \dif x \nonumber\\
   &\qquad\qquad\qquad\qquad\qquad =\int_{\mathbb{T}^3} \int_{\mathbb{R}^3} \left[\frac{|v|^2}2+\Phi(x)\right] \left\{F_0(x,v)-\mu_E(x,v)\right\} \dif v \dif x=E_0,\label{eq:conservation of energy}
  \end{align}
  as well as the entropy inequality:
  \begin{equation}\label{ineq:entropy inequality}
   \mathscr{H}(F(t))-\mathscr{H}(\mu_E) \le \mathscr{H}(F_0)-\mathscr{H}(\mu_E),
  \end{equation}
  where $\mathscr{H}(F)\triangleq\int_{\mathbb{T}^3} \int_{\mathbb{R}^3} F\ln F \dif v \dif x$.

As pointed out in \cite{Kim}, generally, one does not have  the conservation of momentum. In fact, it is easy to derive from the equation that
  \begin{equation}\label{1.10}
   \frac{\dif}{\dif t} \int_{\mathbb{T}^3} \int_{\mathbb{R}^3} v F(t,x,v) \dif v \dif x+\int_{\mathbb{T}^3} \int_{\mathbb{R}^3} \nabla\Phi(x) F(t,x,v) \dif v \dif x=0,
  \end{equation}
  which implies that the momentum is not conserved generally. However, if the external potential $\Phi$ is independent of some $x_i$, we obtain the conservation of momentum for those $v_i$: that is, if we denote the orthogonal complement space of the vector $\nabla\Phi(x)$ by $\nabla\Phi(x)^\perp$, we define $\Lambda$ the degenerate subspace of $\nabla\Phi$ by
  \begin{equation}\nonumber
   \Lambda\triangleq\cap_{x\in\mathbb{T}^3} \nabla\Phi(x)^\perp.
  \end{equation}
  Upon reorienting the coordinates, we assume that $\Lambda$ is spanned by $\{e_1,\cdots,e_{n_0}\}$ as long as $n_0=\dim\Lambda>0$. In other words, $\frac{\partial \Phi}{\partial x_1}=\cdots=\frac{\partial \Phi}{\partial x_{n_0}}\equiv0$, which together with \eqref{1.10} shows the conservation of momentum for degenerate $\{v_1,\cdots,v_{n_0}\}$
  \begin{equation}\label{eq:degenerate conservation of momentum}
   \int_{\mathbb{T}^3} \int_{\mathbb{R}^3} v_i \left\{F(t,x,v)-\mu_E(x,v)\right\} \dif v \dif x=\int_{\mathbb{T}^3} \int_{\mathbb{R}^3} v_i \left\{F_0(x,v)-\mu_E(x,v)\right\} \dif v \dif x = \mathbf{J}_0,
  \end{equation}
  for $i=1,\cdots,n_0$. We point out that the momentum conservation \eqref{eq:degenerate conservation of momentum} is important for us to obtain the decay estimate.

The Boltzmann equation is a fundamental model in the collisional kinetic theory. There have been extensive results on the studies of well-posedness, especially on the global existence theory and large time asymptotic behavior. Among these works, we only mention some known results with external potential. For the classical Boltzmann equation, we refer the interested readers to \cite{BGGL,D-Lion,Guo-03,Guo1,Guo2,IS,KS,LYY,Uk,UA} and the references therein. For the Boltzmann equation with external potential, the local well-posedness was studied in \cite{Asano2, Drange}; the global well-posedness was established in \cite{Duan, Duan-Strain, Duan-Ukai-Yang-Zhao, Li-Yu,Tabata,Ukai-Yang-Zhao, Yu} around a local Maxwellian with small perturbation and some smallness assumption for the external potential $\Phi$ by the energy method. Guo \cite{Guo-04} investigated the global well-posedness of Vlasov-Poisson-Boltzmann solution with a small self-contained external potential, see also \cite{DuanYZ,DuanYZ-1} and the references therein. For the  large amplitude external potential case, it is hard to use the energy method to study the global existence of solution to the Boltzmann equation since the derivatives of Boltzmann solution may grow in time. Recently, Guo \cite{Guo1,Guo2} developed a new $L^2\cap L^\infty$ theory to study the global stability of Boltzmann solutions, which can avoid the derivatives estimates. By using the $L^2\cap L^\infty$ method, Kim \cite{Kim} proved the global well-posedness of Boltzmann equation and large time behavior with large external potential under small perturbation around the local Maxwellian. It is noted that the characteristic lines are heavily bended in \cite{Kim}, then to use the $L^2\cap L^\infty$ method, one has to analyze the possible singular point for the change of variable
\begin{equation}\nonumber
u \rightarrow X(\tau;s,X(s;t,x,v),u).
\end{equation}
This difficulty was overcome in \cite{Kim} by using Asano's almost transversality result in classical dynamicl system \cite{Asano1}. More precisely, Kim proved
\begin{align}\label{1.13}
\mbox{det} \left\{\frac{\partial X(\tau;s,X(s;t,x,v),u)}{\partial u}\right\} \neq 0,
\end{align}
for almost every $(\tau,s,u)\in(0,s)\times(0,t)\times \mathbb{R}^3$ for all $X(\tau;s,X(s;t,x,v),u)\in\mathbb{T}^3$.
We also apply the key estimate \eqref{1.13} in our present paper.

We remark that in those works in the perturbation framework, the initial data are required to have small amplitude around the Maxwellian due to the difficulty coming from the nonlinear term.
Recently, Duan-Huang-Wang-Yang \cite{Duan-Huang-Wang-Yang} developed a new $L^\infty_xL^1\cap L^\infty_{x,v}$ method to obtain the global well-posedness of Boltzamann equation in $\mathbb{T}^3$ or $\mathbb{R}^3$ for a class of initial data with large amplitude oscillations. 
Our aim in the paper is to extend the result of \cite{Kim} to a class of initial data with large amplitude oscillations in $L^\infty_{x,v}$. Due to the appearance of large external potential, even though it is still in periodic box $\mathbb{T}^3$,  one is hard to prove the existence of a global solution to the Boltzmann equation by using the idea of \cite{Duan-Huang-Wang-Yang}. In fact, one of the main ideas in \cite{Duan-Huang-Wang-Yang} is to bound the $L^\infty_x L^1_v$ norm of the solution $f(t,x,v)$ after the local existence time $t_1(\cong\frac1{1+A_0})$ which depends only on the initial $L^\infty$ bound $A_0$, i.e.,
\begin{align}\label{1.12}
\int_{\mathbb{R}^3} |f(t,x,v)| \dif v\le \int_{\mathbb{R}^3}|f_0(x-vt,v)| \dif v+\cdots\le \frac{1}{t_1^2} \|f_0\|_{L^1_x L^\infty_v}+\cdots,\ \forall\quad t\ge t_1,\ x\in\mathbb{T}^3,
\end{align}
where $\|f_0\|_{L^1_xL^\infty_v}$ is demanded to be sufficiently small. Generally, it is difficult to carry the same procedure in \cite{Duan-Huang-Wang-Yang} since the external force could significantly bend the characteristic lines, then one can not use the change of variable as in \eqref{1.12} at all time $t\geq t_1>0$. Recently, Duan-Wang \cite{Duan-Wang} and Duan-Huang-Wang-Zhang \cite{DHWZ} studied the global existence of solution to the Boltzmann equation in bounded domain for a class of initial data with large amplitude oscillations in $L^\infty_{x,v}$, where the similar problem also occurs due to the complex reflections of particles at physical boundary, and they overcame this problem basing on a new Gronwall-type argument. We point out that the perturbation of \cite{DHWZ} is made around a non-constant equilibrium (neither global Maxwellian nor local Maxwellian) which is induced by a given boundary temperature with small oscillations around constant temperature, and it is a big open problem to show the existence and dynamical stability of non-constant equilibriums for large-oscillating boundary temperature. In the current paper, the situation is quite different, as the non-constant equilibrium state induced by the stationary potential of large-amplitude naturally exists, and the perturbation approach developed in \cite{DHWZ,Duan-Wang} is still applicable. Motivated by \cite{Duan-Wang,DHWZ}, we shall use the Gronwall-type argument to investigate the existence of global solution to the Boltzmann equation with large external potential for a class of  initial data with large oscillations.

\

For later use, we define the weight function
\begin{equation}\nonumber
w_\beta(x,v)=\left(\frac{|v|^2}2+\Phi(x)+1\right)^\frac\beta2 \cong \left(\frac{|v|^2}2+1\right)^\frac\beta2.
\end{equation}
Our main result is
  \begin{theorem}\label{thm:main theorem}
   Let $0\le\gamma\le1$, $\beta\ge4$, and $\Phi$ be a periodic $C^3(\mathbb{T}^3)$ function with \eqref{def:the C3 norm of Phi}. Assume $F_0(x,v)=\mu_E(x,v)+\sqrt{\mu_E(x,v)}f_0(x,v)\ge0$, the conservations of mass \eqref{eq:conservation of mass}, energy \eqref{eq:conservation of energy} and momentum for degenerate $\{v_1,\cdots,v_{n_0}\}$ \eqref{eq:degenerate conservation of momentum} are valid for initial data $F_0$ with
   \begin{equation}\nonumber
    (M_0,\mathbf{J}_0,E_0)=(0,0,0) \in \mathbb{R}\times\mathbb{R}^{n_0}\times\mathbb{R}.
   \end{equation}
   For any given $A_0\ge1$, there exists $\kappa_0>0$, depending only on $\beta, M$ and $A_0$, such that if
   \begin{equation}\nonumber
    {\left\|w_\beta f_0 \right\|}_{L^\infty} \le A_0 \quad\text{and}\quad \|f_0\|_{L^2} \le \kappa_0,
   \end{equation}
   the Boltzmann equation \eqref{eq:main} admits a unique global solution $F(t,x,v)=\mu_E(x,v)+\sqrt{\mu_E(x,v)}f(t,x,v)\ge0$ satisfying
   \begin{equation}\nonumber
    \left\|w_\beta f(t)\right\|_{L^\infty} \le \tilde{C}_0 A_0^2\exp\left\{\frac{\tilde{C}_0}{\tilde{\nu}_0}A_0^2\right\} e^{-\tilde\lambda_0 t},
   \end{equation}
   for all $t\ge0$, where $\tilde{C}_0, \tilde{\nu}_0$ and $\tilde\lambda_0$ are some positive constants. In addition, the conservations of mass \eqref{eq:conservation of mass}, energy \eqref{eq:conservation of energy}, and momentum \eqref{eq:degenerate conservation of momentum} as well as the additional entropy inequality \eqref{ineq:entropy inequality} hold.
  \end{theorem}

The paper is organized as follows. In section \ref{section2}, we present some results which will be used later. The proof of main theorem \ref{thm:main theorem} is given in section \ref{sec:3}. We prove the local existence of solution to Boltzmann equation with large external potential in the appendix.\vspace{2mm}

\noindent{\bf Notations:} Throughout the paper, $C$ denotes a positive constant depending only on $\beta,\gamma,M$ which may vary from line to line. $C_{a,b}$ denotes a positive constant if it further depends on $a,b$. For convenience, we denote $\|\cdot\|_{L^p(\mathbb{T}^3\times\mathbb{R}^3)}$ by $\|\cdot\|_{L^p}$ for the rest of the paper.

\

 \section{Preliminaries}\label{section2}
  \subsection{Backward characteristics}
   Given $(t,x,v)\in(0,+\infty)\times\mathbb{T}^3\times\mathbb{R}^3$, let $[X(s),V(s)]$ be the backward bi-characteristics for the Boltzmann equation \eqref{eq:main}, which is determined by the following ODEs
   \begin{equation}\label{eq:characteristics}
    \begin{cases}
     \frac{\dif X}{\dif s}=V(s), \\
     \frac{\dif V}{\dif s}=-\nabla\Phi(X(s)),
     \\ X(t)=x, V(t)=v.
    \end{cases}
   \end{equation}
   It is noted that there exists a Hamiltonian to the system \eqref{eq:characteristics} given by
   \begin{equation}\nonumber
    H(x,v)=\frac{|v|^2}2+\Phi(x).
   \end{equation}
   One can directly check out that the Hamiltonian $H(x,v)$ is preserved along the characteristics, i.e,
   \begin{equation}\nonumber
    \frac{|V(s)|^2}2+\Phi(X(s))=\frac{|v|^2}2+\Phi(x),
   \end{equation}
   which, together with \eqref{def:the C3 norm of Phi}, yields that
   \begin{equation}\label{est:|V(s)|}
    |v|^2-2M \le |V(s)|^2=|v|^2+2\Phi(x)-2\Phi(X(s)) \le |v|^2+2M.
   \end{equation}
   With the above characteristic line, one can express the mild solution of \eqref{eq:equation of f} as
   \begin{equation}\label{eq:mild form of f}
    \begin{split}
     f(t,x,v)&=e^{-\int_0^t g(\tau) \dif \tau} f_0(X(0),V(0))+\int_0^t e^{-\int_s^t g(\tau) \dif \tau-\Phi(X(s))} Kf(s,X(s),V(s)) \dif s
     \\&\quad+\int_0^t e^{-\int_s^t g(\tau) \dif \tau-\frac{\Phi(X(s))}2} \Gamma(f,f)(s,X(s),V(s)) \dif s,
    \end{split}
   \end{equation}
   where $g$ is defined as
   \begin{equation}\label{def:g}
    g(\tau)=e^{-\Phi(X(\tau))}\nu(V(\tau)).
   \end{equation}

   \

   By using the almost transversality theorem of \cite{Asano1}, C. Kim \cite{Kim} proved the following useful lemma for Boltzmann equation with steady large external potential.
   \begin{lemma}[Kim \cite{Kim}]\label{lem:Change of Variable}
    Assume $\Phi$ is a periodic $C^3$-function on $\mathbb{T}^3$. Fix $\varepsilon>0$, $t_0>0$ and $N>0$. There are disjoint open interval partitions of the time interval $[0,t_0]:\mathscr{D}_{i^1}^1 \subset [0,t_0]$ for $i^1 \in \{1,2,\cdots,M_1\}$ and disjoint open box partitions of $[-N,N]^3:\mathscr{D}_{I^2}^2 \subset [-N,N]^3$ for $I^2=(i_1^2,i_2^2,i_3^2) \in \{1,2,\cdots,M_2\}^3$ and disjoint open box partitions of $[-N,N]^3:\mathscr{D}_{I^3}^3 \subset [-N,N]^3$ for $I^3=(i_1^3,i_2^3,i_3^3) \in \{1,2,\cdots,M_3\}^3$. For each $i^1,I^2,I^3$ we have $t_{j,i^1,I^2,I^3}\in\mathscr{D}_{i^1}^1$ for $j=1,2,3$ so that
    \begin{equation}\nonumber
     \left\{s \in \mathscr{D}_{i^1}^1:\det\left(\frac{\partial X}{\partial v}\right)(s;t_0,x,v)=0\right\} \subset \cup_{j=1}^3\left(t_{j,i^1,I^2,I^3}-\frac{\varepsilon}{4M_1},t_{j,i^1,I^2,I^3}+\frac{\varepsilon}{4M_1}\right),
    \end{equation}
    for all $(x,v)\in\mathscr{D}_{I^2}^2 \times \mathscr{D}_{I^3}^3$. Moreover, there exists a positive constant $\delta_*=\delta_*(\varepsilon, M_1,M_2,M_3,N,t_0)>0$ such that
    \begin{equation}\label{2.5}
     \left|\det\left(\frac{\partial X}{\partial v}\right)(s;t_0,x,v)\right|>\delta_*,\quad \forall s \notin \cup_{j=1}^3 \left(t_{j,i^1,I^2,I^3}-\frac{\varepsilon}{4M_1},t_{j,i^1,I^2,I^3}+\frac{\varepsilon}{4M_1}\right),
    \end{equation}
    if $(s,x,v)\in\mathscr{D}_{i^1}^1 \times \mathscr{D}_{I^2}^2 \times \mathscr{D}_{I^3}^3$ for all $i^1,I^2,I^3$.
   \end{lemma}

When the initial data is a small perturbation around $\mu_E$ in weighted $L^\infty_{x,v}$, Kim  \cite{Kim} proved the global existence and large time behavior of Boltzmann solution. For later use, we introduce his main result in the following:
   \begin{theorem}[Kim \cite{Kim}]\label{thm:kim}
    Let $\beta>3$, $\Phi$ be a periodic $C^3$-function on $\mathbb{T}^3$ and $\Phi(x)=\Phi(x_{n_0+1},\cdots,x_3)$ for some $n_0\le3$. Assume that the conservations of mass \eqref{eq:conservation of mass}, energy \eqref{eq:conservation of energy} and momentum for degenerate $\{v_1,\cdots,v_{n_0}\}$ \eqref{eq:degenerate conservation of momentum} are valid for initial data $F_0=\mu_E+\sqrt{\mu_E}f_0$ with
    \begin{equation}\nonumber
     (M_0,\mathbf{J}_0,E_0)=(0,0,0)\in \mathbb{R}\times\mathbb{R}^{n_0}\times\mathbb{R},
    \end{equation}
    then there exist $\lambda_0>0$ and small $\delta>0$ such that if $\|w_\beta f_0\|_{L^\infty}\le\delta$, there exists a unique global solution $F(t,x,v)=\mu_E+\sqrt{\mu_E} f(t,x,v)\ge0$ for the Boltzmann equation \eqref{eq:main} with
    \begin{equation}\nonumber
     \sup_{0\le t<+\infty} \{e^{\lambda_0 t} \|w_\beta f(t)\|_{L^\infty}\}\leq C_0\|w_\beta f_0\|_{L^\infty},
    \end{equation}
    where $C_0$ is some positive constant depending only on $M$.
   \end{theorem}

  \

  \subsection{Useful inequalities}
   Recall $K=K_2-K_1$ in \eqref{def:K1} and \eqref{def:K2}. Let $k(\cdot,\cdot)$ be the kernel of $K$, i.e.,
   \begin{equation}\nonumber
    Kf(v)=\int_{\mathbb{R}^3} k(v,u)f(u) \dif u,
   \end{equation}
   with the symmetric property $k(v,u)=k(u,v)$.
   \begin{lemma}[\cite{Glassey,Guo1}]\label{lem:2.2}
    For $0\le\gamma\le1$, one has
    \begin{equation}\label{est:Grad's estimate}
     |k(v,u)|\le C\left\{|v-u|+|v-u|^{-1}\right\} e^{-\frac{|v-u|^2}8} e^{-\frac{||v|^2-|u|^2|^2}{8|v-u|^2}}.
    \end{equation}
   \end{lemma}

   By the same calculations as in \cite{Guo1}, it is straightforward to check that for $\alpha\ge0$ and $\theta\in[0,\frac18)$, the following inequality holds for some positive constant $C_\alpha$
   \begin{equation}\label{est:k}
    \int_{\mathbb{R}^3} \left|k(v,u)\cdot \frac{(1+|v|)^{\alpha} e^{\theta |v|^2}}{(1+|u|)^{\alpha} e^{\theta |u|^2}}\right| \dif u \le C_\alpha (1+|v|)^{-1}.
   \end{equation}

   \begin{lemma}[Duan-Wang \cite{Duan-Wang}]
    For $0\le\gamma\le1$, there is a generic constant $C_\alpha>0$ such that
    \begin{equation}\label{est:Gamma+}
     \left|(1+|v|)^{\alpha}\Gamma_{+}(f,f)(v)\right| \le \frac{C_\alpha \|(1+|\cdot|)^{\alpha} f\|_{L_{v}^\infty}}{1+|v|} \left(\int_{\mathbb{R}^3} (1+|u|)^{4}|f(u)|^2 \dif u\right)^{\frac12},
    \end{equation}
    for all $v\in \mathbb{R}^3$. In particular, for $\alpha\geq4$, one has
    \begin{equation}\label{est:Gamma+ alpha>=2}
     \left|(1+|v|)^{\alpha}\Gamma_{+}(f,f)(v)\right| \le \frac{C_\alpha\|(1+|\cdot|)^{\alpha}f\|^2_{L^\infty}}{1+|v|},
    \end{equation}
    for all $v\in\mathbb{R}^3$.
   \end{lemma}

   \

   \begin{lemma}[\cite{Duan-Ukai-Yang-Zhao,UY-2006}]\label{lem:L2 estimate for Gamma}
    Assume $0\le\gamma\le1$, then it holds that
    \begin{equation}\label{est:L2 estimate for Gamma}
     \left\|\Gamma(f,f)\right\|_{L_v^2} \le C \left\|f\right\|_{L_v^2} \left\|\nu f\right\|_{L_v^2},
    \end{equation}
    where the positive constant $C$ depends only on $\gamma$.
   \end{lemma}

 \

 \section{Global Stability}\label{sec:3}
  In this section, we consider the global stability of Boltzmann equation \eqref{eq:main}. By using  the local existence of unique solution to \eqref{eq:main} with arbitrarily large initial data established in Theorem \ref{thm:local existence}, it suffices to obtain uniform estimates on solutions.

  Recall that the initial data satisfies $\|w_\beta f_0\|_{L^\infty}\le A_0$ for a given positive constant $A_0\ge 1$ which could be large. Let $f(t,x,v)$ be the solution to the Boltzmann equation \eqref{eq:equation of f} with initial data $F_0=\mu_E+\sqrt{\mu_E} f_0\ge0$ over the time interval $[0,T)$ for $T\in(0,+\infty)$. Throughout this section, we make the {\it a priori} assumption:
  \begin{equation}\label{a priori assumption}
   \sup_{0\le t<T} \|h(t)\|_{L^\infty} \le A_1,
  \end{equation}
  with $h(t,x,v):=w_\beta(x,v) f(t,x,v)$, where $A_1\ge1$ is a large positive constant depending only on $A_0$,  and will be determined later in section \ref{section3.4}.

  \

It follows from \eqref{eq:equation of f} that
  \begin{equation}\label{eq:equation of h}
   \frac{\partial h}{\partial t} +v \cdot \nabla_x h -\nabla\Phi \cdot \nabla_v h + e^{-\Phi}\nu(v) h =e^{-\Phi}K_\beta h +e^{-\frac\Phi2} \Gamma_\beta(h,h),
  \end{equation}
  where the weighted operator $K_\beta,\Gamma_\beta$ are defined by
  \begin{equation}\nonumber
   K_\beta h=w_\beta K\left(\frac{h}{w_\beta}\right) \quad\text{and}\quad \Gamma_\beta(h,h)=w_\beta \Gamma\left(\frac{h}{w_\beta},\frac{h}{w_\beta}\right).
  \end{equation}
  Integrating along the backward trajectory defined in \eqref{eq:characteristics}, one obtains the mild formula for $h$
  \begin{equation}\label{eq:mild form of h}
   \begin{split}
    h(t,x,v)&=e^{-\int_0^t g(\tau) \dif \tau} h_0(X(0),V(0))+\int_0^t e^{-\int_s^t g(\tau) \dif \tau-\Phi(X(s))} K_\beta h(s,X(s),V(s)) \dif s,\nonumber
    \\&\quad+\int_0^t e^{-\int_s^t g(\tau) \dif \tau-\frac{\Phi(X(s))}2} \Gamma_\beta (h,h)(s,X(s),V(s)) \dif s,
   \end{split}
  \end{equation}
  where $g(\tau)$ is defined in \eqref{def:g}.

  \

  \subsection{$L^2$ a priori estimate}
   \begin{lemma}\label{lem:L2 estimate for f}
    Let $f(t,x,v)$ be the solution to Boltzmann equation \eqref{eq:equation of f}, then there exists a generic positive constant $\tilde{C}_1\ge1$ such that
    \begin{equation}\label{est:L2 estimate for f}
     \|f(t)\|_{L^2} \le \tilde{C}_1 \|f_0\|_{L^2} e^{\tilde{C}_1 A_1 t},\quad \forall t\in[0,T].
    \end{equation}
   \end{lemma}
\noindent{\bf Proof.} Applying the standard energy estimate to \eqref{eq:equation of f}, we deduce
  \begin{equation}\nonumber
     \frac{\dif}{\dif t} \|f(t)\|_{L^2}^2+2(e^{-\Phi}f(t),Lf(t)) \le 2\int_{\mathbb{T}^3} \int_{\mathbb{R}^3} e^{-\frac{\Phi}2} f(t,x,v) \Gamma(f,f)(t,x,v) \dif v \dif x.
    \end{equation}
    Notice that $(e^{-\Phi}f(t),Lf(t))\ge0$, then we deduce from Lemma \ref{lem:L2 estimate for Gamma} and the a priori assumption \eqref{a priori assumption} that
    \begin{equation}\nonumber
     \frac{\dif}{\dif t} \|f(t)\|_{L^2}^2 \le C\|\nu f(t)\|_{L^2} \|f(t)\|_{L^2}^2 \le C A_1 \|f(t)\|_{L^2}^2.
    \end{equation}
    Therefore, \eqref{est:L2 estimate for f} follows from the Gronwall's inequality, and the proof of Lemma \ref{lem:L2 estimate for f} is completed. $\hfill\Box$

  \

  \subsection{$L^\infty_xL^1_v$ estimate}
   To consider the Boltzmann equation with a class of large amplitude initial data, motivated by \cite{Duan-Wang}, it is better to rewrite the Boltzmann equation \eqref{eq:equation of h} as
   \begin{equation}\label{eq:equation of h with loss term on the left}
    \begin{split}
     \frac{\partial h}{\partial t} +v \cdot \nabla_x h -\nabla\Phi \cdot \nabla_v h+R(F)h =e^{-\Phi} K_\beta(h)+e^{-\frac\Phi2} w_{\beta} \Gamma_+\left(\frac{\sqrt\mu h}{w_\beta},\frac{\sqrt\mu h}{w_\beta}\right),
    \end{split}
   \end{equation}
   where
   \begin{equation}\label{def:R(F)}
    \begin{split}
     R(F)&:=\int_{\mathbb{R}^3} \int_{\mathbb{S}^2} B(v-u,\omega) F(t,x,u) \dif \omega \dif u
     \\&=\int_{\mathbb{R}^3} \int_{\mathbb{S}^2} B(v-u,\omega) [\mu_{E}(x,u)+\sqrt{\mu_E(x,u)}f(t,x,u)] \dif \omega \dif u.
    \end{split}
   \end{equation}
   It follows from the fact $F(t,x,v)\ge0$ that
   \begin{equation}\nonumber
    I(t,s):=\exp\left\{-\int_s^t R(F)(\tau,X(\tau),V(\tau)) \dif \tau\right\}\le1.
   \end{equation}
   In order to extend the local-in-time solution to a global one with possible large data in $L^\infty_{x,v}$, it is necessary to further obtain the time-decay property of $I(t,s)$, namely,
   \begin{equation}\label{est:I(t,s)}
    I(t,s) \le Ce^{-\frac1C (t-s)},
   \end{equation}
   for some generic large positive constant $C$. For solution of large amplitude, it is impossible to obtain the above time-decay property initially since the vacuum could not be excluded. However, if the $L^2$ norm of initial data is small, such time-decay property may still hold in some sense even if the initial data is allowed to have large oscillations.

   To prove \eqref{est:I(t,s)}, one needs to recover a positive lower bound for $R(F(t))$. Indeed, we aim to obtain a positive lower bound for $R(F(t))$ for time $t$ suitably large. In fact, we have the following lemma.
   \begin{lemma}\label{lem:3.2}
    Under the a priori assumption \eqref{a priori assumption}, there exists a constant $\tilde{C}_2\ge1$ such that for any given positive time $T_1>\tilde{t}$ with
    \begin{equation}\label{3.7}
     \tilde{t}:=\frac{2e^M}{\nu_0} \ln(\tilde{C}_2e^{\frac{M}{2}}A_0),
    \end{equation}
    there is a small positive constant $\kappa_1=\kappa_1(A_1,T_1)>0$, depending only on $A_1$ and $T_1$, such that if $\|f_0\|_{L^2}\le\kappa_1$, one has
    \begin{equation}\label{3.8}
     R(F(t))\ge\frac12 e^{-\Phi} \nu(v),
    \end{equation}
    for all $t\in[\tilde{t},\min\{T,T_1\}]$, where $\kappa_1$ decreases in $A_1$ and $T_1$.
   \end{lemma}
\noindent{\bf Proof.} To obtain \eqref{3.8}, it is noted that
    \begin{equation}\nonumber
     \int_{\mathbb{R}^3} \int_{\mathbb{S}^2} B(v-u,\omega) \sqrt{\mu(u)} |f(t,x,u)| \dif u \le C_1 \nu(v) \int_{\mathbb{R}^3} e^{-\frac{|u|^2}8} |f(t,x,u)| \dif u,
    \end{equation}
    where $C_1\ge1$ is some generic constant. Hence, by using \eqref{def:R(F)}, it suffices to prove
    \begin{equation}\label{3.9}
     \int_{\mathbb{R}^3} e^{-\frac{|v|^2}8} |f(t,x,v)| \dif v \le \frac1{2C_1} e^{-\frac M2}.
    \end{equation}
    It follows from \eqref{eq:mild form of f} that
    \begin{equation}\label{3.10}
     \begin{split}
      &\int_{\mathbb{R}^3} e^{-\frac{|v|^2}8} |f(t,x,v)| \dif v
      \\&\le \int_{\mathbb{R}^3} e^{-\frac{|v|^2}8} e^{-\int_0^t g(\tau) \dif \tau} |f_0(X(0),V(0))| \dif v
      \\& \quad + \int_{\mathbb{R}^3} e^{-\frac{|v|^2}8} \int_0^t e^{-\int_s^t g(\tau) \dif \tau-\Phi(X(s))} |Kf(s,X(s),V(s))| \dif s \dif v
      \\& \quad + \int_{\mathbb{R}^3} e^{-\frac{|v|^2}8} \int_0^t e^{-\int_s^t g(\tau) \dif \tau-\frac{\Phi(X(s))}2} |\Gamma(f,f)(s,X(s),V(s))| \dif s \dif v
      \\&:=I_1+I_2+I_3.
     \end{split}
    \end{equation}
    Recall $\nu_0=\inf_{v\in\mathbb{R}^3}\nu(v)>0$ and note that
    \begin{equation}\nonumber
     e^{-\int_s^t g(\tau) \dif \tau} \le \exp\left\{-e^{-M}\nu_0(t-s)\right\},
    \end{equation}
it is easy to obtain that
    \begin{equation}\label{est:I1}
     I_1 \le C \|f_0\|_{L^\infty} \exp\left\{-e^{-M} \nu_0 t\right\}.
    \end{equation}

\

    For $I_2$, we split it into three parts.
    \begin{equation}\label{ineq:I2}
     \begin{split}
      I_2 & \le \int_0^t \int_{|v|\ge N} e^{-\frac{|v|^2}8} \exp\left\{-e^{-M}\nu_0(t-s)\right\} |Kf(s,X(s),V(s))| \dif v \dif s
      \\&\quad+\int_0^t \int_{|v|<N} e^{-\frac{|v|^2}8} \exp\left\{-e^{-M}\nu_0(t-s)\right\} \int_{|u|\ge 3N} |k(V(s),u)f(s,X(s),u)| \dif u \dif v \dif s
      \\&\quad+\int_0^t \int_{|v|<N} e^{-\frac{|v|^2}8} \exp\left\{-e^{-M}\nu_0(t-s)\right\} \int_{|u|<3N} |k(V(s),u)f(s,X(s),u)| \dif u \dif v \dif s
      \\&:=I_{21}+I_{22}+I_{23}.
     \end{split}
    \end{equation}
    Using \eqref{est:k}, one obtain that
    \begin{equation}\label{est:I21}
     I_{21} \le C e^{-\frac{N^2}{16}} \sup_{0\le s\le t} \|f(s)\|_{L^\infty}.
    \end{equation}
    It follows from \eqref{est:|V(s)|} that
    \begin{equation}\nonumber
     |V(s)-u| \ge |u|-|V(s)| \ge N,\quad \text{if} \quad |v|\le N, |u|\ge3N, N\gg M,
    \end{equation}
    which yields that 	
    \begin{equation}\label{est:I22}
     \begin{split}
      I_{22} &\le Ce^{-\frac{N^2}{32}} \sup_{0\le s\le t} \|f(s)\|_{L^\infty} \int_0^t \int_{|v|<N} e^{-\frac{|v|^2}8} \exp\left\{-e^{-M}\nu_0(t-s)\right\}
      \\&\qquad\qquad\qquad\qquad\qquad\qquad\times\int_{|u|\ge3N} |k(V(s),u)|e^{\frac{|V(s)-u|^2}{32}} \dif u \dif v \dif s
      \\&\le Ce^{-\frac{N^2}{32}} \sup_{0\le s\le t} \|f(s)\|_{L^\infty},
     \end{split}
    \end{equation}
    where we have used the following fact
    \begin{equation}\nonumber
     \int_{\mathbb{R}^3} |k(V(s),u)|e^{\frac{|V(s)-u|^2}{32}} \dif u \leq \frac{C}{1+|V(s)|} \cong \frac{C}{1+|v|}.
    \end{equation}
    To estimate $I_{23}$, it follows from \eqref{est:Grad's estimate} and the H\"{o}lder's inequality that
    \begin{equation}\label{3.15}
     \begin{split}
      I_{23} &\le C \int_0^t \exp\{-e^{-M} \nu_0 (t-s)\}\bigg\{\int_{|v|<N} \int_{|u|<3N} e^{-\frac{|v|^2}8-\frac{|V(s)-u|^2}8}
      \\&\qquad\qquad\qquad\qquad\times|f(s,X(T_1-t+s;T_1,x,v),u)|^2 \dif u \dif v\bigg\}^\frac12 \dif s,
     \end{split}
    \end{equation}
    where we have used the fact
    \begin{equation}\label{eq:translation invariance of characteristics}
     X(s)=X(s;t,x,v)=X(T_1-t+s;T_1,x,v),
    \end{equation}
    since the external potential $\Phi$ is time-independent. We apply Lemma \ref{lem:Change of Variable} for the case $t_0=T_1$ to deal with the term on the right hand side of \eqref{3.15}. Assume $x\in\mathscr{D}_{I_2}^2$ for some $I_2\in\{1,2,\cdots,M_2\}^3$, then it follows from  Lemma \ref{lem:Change of Variable} that
    \begin{align}\label{ineq:I23}
      I_{23} &\le C \sum_{i_1=1}^{M_1} \sum_{I_3\in\{1,\cdots,M_3\}^3} \int_0^t \exp\{-e^{-M} \nu_0 (t-s)\} \mathbf{1}_{\mathscr{D}_{i^1}^1}(T_1-t+s)\nonumber\\
      &\qquad\qquad\qquad\qquad\qquad\times\left\{\int_{v\in\mathscr{D}_{I^3}^3} \int_{|u|\le3N} \cdots \dif v \dif u\right\}^{\frac12}\dif s\nonumber\\
      &=C\sum_{i^1} \sum_{j=1}^3 \sum_{I^3} \int_0^t \mathbf{1}_{\mathscr{D}_{i^1}^1\cap(t_{j,i_1,I_2,I_3}-\frac{\varepsilon}{4M_1},t_{j,i_1,I_2,I_3}+\frac{\varepsilon}{4M_1})}(T_1-t+s) \nonumber\\
      &\qquad\qquad\qquad\times \exp\{-e^{-M} \nu_0 (t-s)\}\cdot\left\{\int_{|v|\le N} \int_{|u|\le3N} \mathbf{1}_{\mathscr{D}_{I_3}^3}(v)\cdots \dif v \dif u\right\}^{\frac12}\dif s \nonumber\\
      &\quad+C\sum_{i^1} \sum_{j=1}^3 \sum_{I_3} \int_0^t  \mathbf{1}_{\mathscr{D}_{i^1}^1\setminus (t_{j,i_1,I_2,I_3}-\frac{\varepsilon}{4M_1},t_{j,i_1,I_2,I_3}+\frac{\varepsilon}{4M_1})}(T_1-t+s)
      \nonumber\\
      &\qquad\qquad\qquad \times \exp\{-e^{-M} \nu_0 (t-s)\}\cdot\left\{\int_{|v|\le N} \int_{|u|\le3N} \mathbf{1}_{\mathscr{D}_{I_3}^3}(v)\cdots \dif v \dif u\right\}^{\frac12}\dif s.
    \end{align}
    A direct calculation shows that the first term on the RHS of \eqref{ineq:I23} is bounded by
    \begin{equation}\label{est:the first term of I23}
     C\sum_{i_1=1}^{M_1}\sum_{j=1}^3 \frac{\v}{2M_1} \sup_{0\le s\le t} \|f(s)\|_{L^\infty} \le C  \varepsilon \sup_{0\le s\le t} \|f(s)\|_{L^\infty}.
    \end{equation}
    To estimate the second term on the RHS of \eqref{ineq:I23}, we consider the following change of variables
    \begin{equation}\label{3.19}
     v \rightarrow X(T_1-t+s;T_1,x,v).
    \end{equation}
    It follows from Lemma \ref{lem:Change of Variable} that
    \begin{equation}\label{3.20}
     \left|\frac{\partial X(T_1-t+s; T_1,x,v)}{\partial v}\right| \ge \delta_*(\v, M_1,M_2,M_3,N,T_1)>0,
    \end{equation}
    for $T_1-t+s\in\mathscr{D}_{i^1}^1\setminus (t_{j,i_1,I_2,I_3}-\frac{\varepsilon}{4M_1},t_{j,i^1,I^2,I^3}+\frac{\varepsilon}{4M_1})$ with $j=1,2,3$ and $i_1=1,\cdots, M_1$. Applying the change of variables \eqref{3.19} and using \eqref{3.20}, the second term on the RHS of \eqref{ineq:I23} is bounded as
    \begin{equation}\label{est:the second term of I23}
     C(\varepsilon,M_1,M_2,M_3,N,T_1) \int_0^t \exp\{-e^{-M} \nu_0 (t-s)\} \|f(s)\|_{L^2} \dif s,
    \end{equation}
    which, together with \eqref{est:the first term of I23}, \eqref{ineq:I23}, \eqref{est:I22}, \eqref{est:I21} and \eqref{ineq:I2}, yields that
    \begin{equation}\label{est:I2}
     \begin{split}
      I_2&\le C\left(\varepsilon+e^{-\frac{N^2}{32}}\right) \sup_{0\le s\le t} \|f(s)\|_{L^\infty}
      \\&\quad+C(\varepsilon,M_1,M_2,M_3,N,T_1) \int_0^t \exp\{-e^{-M} \nu_0 (t-s)\} \|f(s)\|_{L^2} \dif s.
     \end{split}
    \end{equation}

 For $I_3$, we split it into several parts
    \begin{equation}\label{ineq:I3}
     \begin{split}
      I_3 &\le \int_0^t \int_{|v|\ge N} e^{-\int_s^t g(\tau) \dif \tau-\frac{\Phi(X(s))}2} e^{-\frac{|v|^2}8} |\Gamma(f,f)(s,X(s),V(s))| \dif v \dif s
      \\& \quad +\int_0^t \int_{|v|<N} e^{-\int_s^t g(\tau) \dif \tau-\frac{\Phi(X(s))}2} e^{-\frac{|v|^2}8} \int_{|u|\ge3N} (\Lambda_++\Lambda_-)(s,t,x,v,u) \dif u \dif v \dif s
      \\&\quad +\int_0^t \int_{|v|<N} \exp\{-e^{-M} \nu_0 (t-s)\} e^{-\frac{|v|^2}8} \int_{|u|<3N} \Lambda_-(s,t,x,v,u) \dif u \dif v \dif s
      \\&\quad +\int_0^t \int_{|v|<N} \exp\{-e^{-M} \nu_0 (t-s)\} e^{-\frac{|v|^2}8} \int_{|u|<3N} \Lambda_+(s,t,x,v,u) \dif u \dif v \dif s
      \\&:=I_{31}+I_{32}+I_{33}+I_{34},
     \end{split}
    \end{equation}
    where
    \begin{equation}\nonumber
     \begin{split}
     &\Lambda_+(s,t,x,v,u)=\int_{\mathbb{S}^2} B(V(s)-u,\omega) e^{-\frac{|u|^2}4} |f(s,X(s),u')f(s,X(s),v')| \dif \omega,
     \\&\Lambda_-(s,t,x,v,u)=\int_{\mathbb{S}^2} B(V(s)-u,\omega) e^{-\frac{|u|^2}4} |f(s,X(s),u)f(s,X(s),V(s)) \dif \omega,
    \end{split}
   \end{equation}
   with $u'=u-[(u-V(s))\cdot\omega]\omega$ and $v'=V(s)+[(u-V(s))\cdot \omega]\omega$.\vspace{3mm}

Noting the fact $\nu(V(s))\cong \nu(v)$, we have that
   \begin{equation} \label{3.23-1}
    |\Gamma(f,f)(s,X(s),V(s))|\le C\nu(V(s)) \|f(s)\|_{L^\infty}^2\leq C\nu(v) \|f(s)\|_{L^\infty}^2,
   \end{equation}
   and
   \begin{equation} \label{3.23-2}
    \int_0^t e^{-\int_s^t g(\tau) \dif \tau-\frac{\Phi(X(s))}2} \nu(V(s)) \dif s \leq C.
   \end{equation}
Now by using \eqref{3.23-1} and \eqref{3.23-2}, it holds that
   \begin{equation}\label{est:I31}
    \begin{split}
     I_{31} &\le C \sup_{0\le s\le t} \|f(s)\|_{L^\infty}^2 \int_{|v|\ge N} e^{-\frac{|v|^2}8} \dif v \int_0^t e^{-\int_s^t g(\tau) \dif \tau-\frac{\Phi(X(s))}2} \nu(V(s)) \dif s
     \\&\le C e^{-\frac{N^2}{16}} \sup_{0\le s\le t} \|f(s)\|_{L^\infty}^2,
    \end{split}
   \end{equation}
   and
   \begin{equation}\label{est:I32}
    \begin{split}
     I_{32} &\le C \sup_{0\le s\le t} \|f(s)\|_{L^\infty}^2 \int_0^t \int_{|v|<N} e^{-\int_s^t g(\tau) \dif \tau-\frac{\Phi(X(s))}2} e^{-\frac{|v|^2}8} \int_{|u|\ge3N} |V(s)-u|^\gamma e^{-\frac{|u|^2}4} \dif u \dif v \dif s
     \\& \le C e^{-\frac{9N^2}8} \sup_{0\le s\le t} \|f(s)\|_{L^\infty}^2 \int_{|v|<N} e^{-\frac{|v|^2}8} \int_0^t e^{-\int_s^t g(\tau) \dif \tau-\frac{\Phi(X(s))}2} \nu(V(s)) \dif s\dif v
     \\&\le C e^{-\frac{9N^2}8} \sup_{0\le s\le t} \|f(s)\|_{L^\infty}^2.
    \end{split}
   \end{equation}
   For $I_{33}$, it follows from the Holder's inequality and similar arguments as in \eqref{3.15}-\eqref{est:the second term of I23} that
   \begin{equation}\label{est:I33}
    \begin{split}
     I_{33} &\le C\sup_{0\le s\le t} \|f(s)\|_{L^\infty}  \int_0^t \exp\{-e^{-M} \nu_0 (t-s)\}
     \\&\qquad\qquad\quad\times\left(\int_{|v|<N} \int_{|u|<3N} e^{-\frac{|v|^2}8} e^{-\frac{|u|^2}4} |f(s,X(T_1-t+s;T_1,x,v),u)|^2 \dif u \dif v\right)^\frac12
     \\&\le C\sup_{0\le s\le t} \|f(s)\|_{L^\infty}  \Bigg\{\varepsilon \sup_{0\le s\le t} \|f(s)\|_{L^\infty}
     \\&\qquad\qquad\qquad\quad+C(\varepsilon,M_1,M_2,M_3,N,T_1) \int_0^t \exp\{-e^{-M} \nu_0 (t-s)\} \|f(s)\|_{L^2} \dif s \Bigg\}
     \\&\le C\varepsilon \sup_{0\le s\le t} \|f(s)\|^2_{L^\infty}
     \\&\qquad\qquad\quad+C(\v,M_1,M_2,M_3,N,T_1) \int_0^t \exp\{-e^{-M} \nu_0 (t-s)\} \|f(s)\|^2_{L^2} \dif s.
    \end{split}
   \end{equation}
   For $I_{34}$, using the standard change of variables, one has
   \begin{equation}\nonumber
    e^{-\frac{|v|^2}{16}} \int_{|u|<3N} \Lambda_+(s,t,x,v,u) \dif u \le C \|f(s)\|_{L^\infty} \left(\int_{|\eta|\le4N} |f(s,X(s),\eta)|^2 \dif \eta\right)^{\frac12},
   \end{equation}
   which, together with similar arguments as in \eqref{3.15}-\eqref{est:the second term of I23}, yields that
   \begin{equation}\label{est:I34}
    \begin{split}
     I_{34} &\le C\left(\varepsilon+e^{-\frac{N^2}{32}}\right) \sup_{0\le s\le t} \|f(s)\|^2_{L^\infty}
     \\&\qquad\qquad\quad+C(\varepsilon,M_1,M_2,M_3,N,T_1) \int_0^t \exp\{-e^{-M} \nu_0 (t-s)\} \|f(s)\|^2_{L^2} \dif s.
    \end{split}
   \end{equation}
   Plug \eqref{est:I31},\eqref{est:I32},\eqref{est:I33} and \eqref{est:I34} into \eqref{ineq:I3} to deduce
   \begin{equation}\label{est:I3}
    \begin{split}
     I_3 &\le C\left(\varepsilon+e^{-\frac{N^2}{32}}\right) \sup_{0\le s\le t} \|f(s)\|^2_{L^\infty}
     \\&\qquad\quad+C(\varepsilon,M_1,M_2,M_3,N,T_1) \int_0^t \exp\{-e^{-M} \nu_0 (t-s)\} \|f(s)\|^2_{L^2} \dif s.
     \end{split}
    \end{equation}
    Substituting \eqref{est:I1}, \eqref{est:I2} and \eqref{est:I3} into \eqref{3.10}, we obtain
    \begin{equation}\label{3.29}
     \begin{split}
      &\int_{\mathbb{R}^3} e^{-\frac{|v|^2}8} |f(t,x,v)| \dif v
      \\& \le C_2 \|f_0\|_{L^\infty} \exp\left\{-e^{-M} \nu_0 t\right\} + C_2 \left(\varepsilon+e^{-\frac{N^2}{32}}\right) \left\{\sup_{0\le s\le t} \|f(s)\|_{L^\infty} +\sup_{0\le s\le t}\|f(s)\|_{L^\infty}^2\right\}
      \\&\quad+C(\varepsilon,M_1,M_2,M_3,N,T_1) \int_0^t \exp\{-e^{-M} \nu_0 (t-s)\} \left[\|f(s)\|_{L^2}+\|f(s)\|^2_{L^2}\right] \dif s,
     \end{split}
    \end{equation}
    where $C_2\ge1$ is a generic positive constant. Take
    \begin{equation}\nonumber
     \tilde{t}:=\frac{2e^M}{\nu_0} \ln(4C_1C_2e^{\frac{M}{2}}A_0),
    \end{equation}
    then the first term of \eqref{3.29} is bounded by
    \begin{equation}\label{3.30}
     C_2 \|f_0\|_{L^\infty} \exp\left\{-e^{-M} \nu_0 t\right\} \le \frac1{4C_1} e^{-\frac M2},\quad \forall\  t\geq\tilde{t}.
    \end{equation}
    On the other hand, by using \eqref{a priori assumption} and \eqref{est:L2 estimate for f}, the remaining term of \eqref{3.29} is bounded as
    \begin{equation}\label{3.31}
     C_2\left(\varepsilon+e^{-\frac{N^2}{32}}\right) \left\{A_1 +A_1^2\right\} +C(\varepsilon,M_1,M_2,M_3,N,T_1) e^{2\tilde{C}_1 A_1 T_1} \left\{\|f_0\|_{L^2}+\|f_0\|^2_{L^2}\right\},
    \end{equation}
    for $t\in[0,T_1]$. We first take $\varepsilon>0$ small enough, then $N\ge1$ suitably large, and finally take initial data $f_0$ such that  $\|f_0\|_{L^2}\le\kappa_1$ with $\kappa_1=\kappa_1(A_1,T_1)>0$ sufficiently small, so that \eqref{3.31} is bounded by
    \begin{equation}\label{3.32}
     C_2\left(\varepsilon+e^{-\frac{N^2}{32}}\right) \left\{A_1 +A_1^2\right\}+C(\v,M_1,M_2,M_3,N,T_1) e^{2\tilde{C}_1A_1T_1} \left\{\kappa_1+\kappa_1^2\right\} \le \frac1{4C_1} e^{-\frac M2}.
    \end{equation}
    Notice that $\kappa_1(A_1,T_1)$ can be decreasing in $A_1$ and $T_1$. Combining \eqref{3.29},  \eqref{3.30}, \eqref{3.31} and \eqref{3.32}, one obtains that
    \begin{equation}\nonumber
     \int_{\mathbb{R}^3} e^{-\frac{|v|^2}8} |f(t,x,v)| \dif v \le \frac1{2C_1} e^{-\frac M2},
    \end{equation}
    which yields \eqref{3.9}. Set $\tilde{C}_2:=4C_1C_2$, then the proof of Lemma \ref{lem:3.2} is completed. $\hfill\Box$

  \

  \subsection{$L^\infty$-estimate}
   In this subsection, we focus on the $L^\infty$ {\it a priori} estimate for $h$. Using \eqref{eq:equation of h with loss term on the left}, we rewrite the mild form of Boltzmann equation as
   \begin{equation}\label{eq:mild form of h with loss term on the left}
    \begin{split}
     &h(t,x,v)=e^{-\int_0^t \tilde{g}(\tau) \dif \tau} h_0(X(0),V(0))
     \\& \quad+\int_0^t \exp\left\{-\int_s^t \tilde{g}(\tau) \dif \tau-\Phi(X(s))\right\} K_\beta h(s,X(s),V(s)) \dif s
     \\&\quad+\int_0^t \exp\left\{-\int_s^t \tilde{g}(\tau) \dif \tau-\frac{\Phi(X(s))}2\right\} w_\beta\Gamma_+\left(\frac{h}{w_\beta},\frac{h}{w_\beta}\right)(s,X(s),V(s)) \dif s,
    \end{split}
   \end{equation}
where
   \begin{equation}\nonumber
    \tilde{g}(\tau)=R(F)(\tau,X(\tau),V(\tau)).
   \end{equation}
To treat the nonlinear $L^\infty$ estimate for the solution $f(t,x,v)$ of Boltzmann equation \eqref{eq:equation of f} with $L^\infty$ large amplitude initial data, one needs to use the time-decay property, i.e., $\tilde{g}(\tau)$ has positive lower bound. Even though $\tilde{g}(0)$ may not have positive lower bound for all time, by using Lemma \ref{lem:3.2}, one obtains that
   \begin{align}\label{dec:exponential decay}
     \exp\left\{-\int_s^t \tilde{g}(\tau) \dif \tau\right\} &\le
     \begin{cases}
      1,& 0\le s\le t\le \tilde{t}, \vspace{1mm}\\
      \exp\left\{-\frac12e^{-M}\nu_0(t-\tilde{t})\right\}, & 0\le s\le \tilde{t}\le t\le T_1, \vspace{1mm}\\
      \exp\left\{-\frac12e^{-M}\nu_0(t-s)\right\}, & 0\le\tilde{t}\le s\le t\le T_1,
     \end{cases}\nonumber\\
     &\le \exp\left\{-\frac12e^{-M}\nu_0(t-s)\right\} \exp\left\{\frac12 e^{-M}\nu_0\tilde{t} \right\}\nonumber\\
     &=\tilde{C}_2 e^{\frac{M}{2}} A_0 e^{-\frac12\tilde{\nu}_{0}(t-s)}, \ 0\leq s\leq t\leq T_1,
   \end{align}
   where $\tilde{t}$ is defined in \eqref{3.7} and we denote $\tilde{\nu}_{0}:=e^{-M}\nu_0$ in the last inequality for simplicity of presentation.
   \begin{lemma}\label{lem:3.3}
    Assume $\|f_0\|_{L^2}\le\kappa_1=\kappa_1(A_1,T_1)$, then there exists a generic positive constant $\tilde{C}_3\ge1$ such that
    \begin{equation}\nonumber
     \begin{split}
      \|h(t)\|_{L^\infty}&\le\tilde{C}_3 A_0^2\left[1+\int_0^t \|h(s)\|_{L^\infty} \dif s \right] \exp\left\{-\frac14 \tilde{\nu}_0 t \right\}
      \\&\quad+A_0 \left(\frac{\tilde{C}_3}{\sqrt{N}}+C_N\cdot\varepsilon\right) \sup_{0\le s\le t}\left\{\|h(s)\|_{L^\infty}+\|h(s)\|_{L^\infty}^3\right\}
      \\&\quad+A_0 C(\varepsilon,M_1,M_2,M_3,N,T_1) \sup_{0\le s\le t} \left\{  \|f(s)\|_{L^2}+\|f(s)\|^3_{L^2}\right\},
     \end{split}
    \end{equation}
    for all $t\in[0,\min\{T,T_1\}]$, where $\varepsilon>0$ is some small parameter to be chosen later, and $N>0$ is some large number determined later.
   \end{lemma}
   \begin{pf}
    It follows from \eqref{eq:mild form of h with loss term on the left} and \eqref{dec:exponential decay} that
    \begin{equation}\label{3.35}
     \begin{split}
      |h(t,x,v)|&\le CA_0^2 e^{-\frac12\tilde{\nu}_{0}t}+C A_0 \int_0^t e^{-\frac12\tilde{\nu}_{0}(t-s)} \int_{\mathbb{R}^3} \left|k_\beta(V(s),u) h(s,X(s),u)\right| \dif u \dif s
      \\&\quad+CA_0 \int_0^t e^{-\frac12\tilde{\nu}_{0}(t-s)} \left|w_\beta \Gamma_+\left(\frac{h}{w_\beta},\frac{h}{w_\beta}\right)(s,X(s),V(s))\right| \dif s
      \\&=CA_0 \left(A_0 e^{-\frac12\tilde{\nu}_{0}t} +J_1+J_2\right),
     \end{split}
    \end{equation}
where $k_{\beta}(v,u)$ is bounded as
    \begin{equation}\nonumber
     |k_\beta(v,u)|\le C\left|k(v,u)\right| \frac{(1+|v|^2)^{\frac{\beta}{2}}}{(1+|u|^2)^{\frac{\beta}{2}}}.
    \end{equation}
    \vspace{1.5mm}

    Now we estimate the terms on the RHS of \eqref{3.35}. For $J_1$, similar as in \cite{Guo2}, we use \eqref{3.35} again to get that
    \begin{equation}\label{ineq:J1}
     \begin{split}
      J_1&\le CA_0^2 \int_0^t e^{-\frac12\tilde{\nu}_{0}(t-s)} \int_{\mathbb{R}^3} \left|k_{\beta}(V(s), u)\right| e^{-\frac12\tilde{\nu}_{0}s} \dif u \dif s
      \\&\quad +CA_0 \int_0^t \int_0^s e^{-\frac12\tilde{\nu}_{0}(t-\tau)} \int_{\mathbb{R}^3}  \int_{\mathbb{R}^3} \left|k_{\beta}(V(s),u) k_{\beta}(\tilde{V}(\tau),\xi) h(\tau,\tilde{X}(s),\xi)\right| \dif \xi \dif u \dif \tau \dif s
      \\&\quad+CA_0 \int_0^t \int_0^s e^{-\frac12\tilde{\nu}_{0}(t-\tau)} \int_{\mathbb{R}^3} \left|k_{\beta}(V(s),u)\right|\cdot \left|w_\beta \Gamma_+\left(\frac{h}{w_\beta},\frac{h}{w_\beta}\right)(\tau,\tilde{X}(\tau),\tilde{V}(\tau))\right|  \dif u \dif \tau \dif s
      \\&:= CA_0 \left(A_0 t e^{-\frac12\tilde{\nu}_{0}t}+J_{11}+J_{12}\right),
     \end{split}
    \end{equation}
    where we have denoted $[\tilde{X}(\tau),\tilde{V}(\tau)]:=[X(\tau;s,X(s),u),V(\tau;s,X(s),u)]$.
    \vspace{1.5mm}

    For $J_{11}$, we divide it into the following several cases.\\
    {\rm Case 1.} For $|v|\ge N\gg1$, it follows from \eqref{est:|V(s)|} and \eqref{est:k} that
    \begin{equation}\label{3.37}
     \begin{split}
      J_{11} &\le \sup_{0\le s\le t} \|h(s)\|_{L^\infty} \int_0^t \int_0^s e^{-\frac12 \tilde{\nu}_0(t-\tau)} \int_{\mathbb{R}^3} \int_{\mathbb{R}^3} \left|k_{\beta}(V(s), u) k_{\beta}(\tilde{V}(\tau), \xi) \right| \dif \xi \dif u \dif \tau \dif s
      \\&\le \sup_{0\le s\le t} \|h(s)\|_{L^\infty} \cdot \sup_{0\le s\le t} \frac{C}{1+|V(s)|}\le \frac{C}{1+|v|} \sup_{0\le s\le t} \|h(s)\|_{L^\infty}
      \\&\le \frac{C}{N}  \sup_{0\le s\le t} \|h(s)\|_{L^\infty}.
     \end{split}
    \end{equation}
    {\rm Case 2.} For either $|v|<N, |u|\ge3N$ or $|u|<3N, |\xi|\ge5N$, it follows from \eqref{est:|V(s)|} that we have either $|V(s)-u|\ge N,\ \forall s\in[0,t]$ or $|\tilde{V}(\tau)-\xi|\ge N,\ \forall \tau\in[0,s]$. Therefore, either of the following is valid
    \begin{equation}\label{3.38}
     \begin{split}
      |k_{\beta}(V(s),u)|\le Ce^{-\frac{N^2}{32}} \left|k_{\beta}(V(s),u) e^{-\frac{|V(s)-u|^2}{32}}\right|,\quad \forall s\in [0,t],
      \\|k_{\beta}(\tilde{V}(\tau),\xi)|\le C e^{-\frac{N^2}{32}} \left|k_{\beta}(\tilde{V}(\tau), \xi) e^{-\frac{|\tilde{V}(\tau)-\xi|^2}{32}}\right|, \quad \forall \tau\in[0,s].
     \end{split}
    \end{equation}
    Hence, by using \eqref{est:Grad's estimate}, a direct calculation shows that
    \begin{equation}\label{3.39}
     \begin{split}
      \int_{\mathbb{R}^3} \left|k_{\beta}(V(s), u) e^{-\frac{|V(s)-u|^2}{32}}\right| \dif u\le C(1+|V(s)|)^{-1}\le C(1+|v|)^{-1},\quad \forall s\in [0,t],
      \\ \int_{\mathbb{R}^3}\left|k_{\beta}(\tilde{V}(\tau), \xi) e^{-\frac{|\tilde{V}(\tau)-\xi|^2}{32}}\right| \dif \xi \le C(1+|\tilde{V}(\tau)|)^{-1}\le C(1+|u|)^{-1}, \quad \forall \tau\in[0,s].
     \end{split}
    \end{equation}
Now it follows from \eqref{3.38} and \eqref{3.39} that
    \begin{equation}\label{3.40}
     \begin{split}
      &\int_0^t \int_0^s e^{-\frac12\tilde{\nu}_{0}(t-\tau)} \left\{\int_{|u|\ge3N} \int_{\mathbb{R}^3}+\int_{|u|<3N} \int_{|\xi|\ge5N}\right\}
      \\&\qquad\qquad\qquad\qquad\times\left|k_{\beta}(V(s),u)k_{\beta}(\tilde{V}(\tau),\xi) h(\tau,\tilde{X}(s),\xi) \right| \dif \xi \dif u \dif \tau \dif s
      \\&\le C e^{-\frac{N^2}{32}} \sup_{0\le s\le t} \|h(s)\|_{L^\infty}.
     \end{split}
    \end{equation}
    {\rm Case 3.} $|v|<N,|u|<3N,|\xi|<5N$. This is the last remaining case. Recalling that $\Phi$ is independent of time, similar as in \eqref{eq:translation invariance of characteristics}, it holds that
    \begin{equation}\label{3.41}
     \tilde{X}(\tau)=X(\tau;s,X(s),u)=X(T_1-s+\tau; T_1,X(s),u).
    \end{equation}
    Using \eqref{est:Grad's estimate}, it is easy to get that
    \begin{equation}\label{3.42}
     \int_{\mathbb{R}^3} \left|k_\beta(v,u)\right|^2 \dif u \le C\int_{\mathbb{R}^3} \left|k(v,u)\right|^2 \frac{(1+|v|^2)^\beta}{(1+|u|^2)^\beta} \dif u \le C(1+|v|)^{-1},
    \end{equation}
    which, together with H\"{o}lder's inequality and \eqref{3.41} yields that
    \begin{equation}\label{3.43}
     \begin{split}
      &\int_0^t\int_0^s e^{-\frac12\tilde{\nu}_{0}(t-\tau)} \int_{|u|<3N} \int_{|\xi|<5N} \left|k_{\beta}(V(s),u)k_{\beta}(\tilde{V}(\tau),\xi) h(\tau,\tilde{X}(s),\xi)\right| \dif \xi \dif u \dif \tau \dif s
      \\&\le \int_0^t \int_0^s e^{-\frac12\tilde{\nu}_{0}(t-\tau)} \left\{\int_{|u|<3N} \int_{|\xi|<5N} \left|k_{\beta}(V(s),u)k_{\beta}(\tilde{V}(\tau),\xi) \right|^2 \dif \xi \dif u \right\}^{\frac12}
      \\&\qquad\qquad\qquad\qquad\qquad\qquad\times \left\{\int_{|u|<3N} \int_{|\xi|<5N} |h(\tau,\tilde{X}(\tau),\xi)|^2 \dif \xi \dif u\right\}^{\frac12} \dif \tau \dif s
      \\&\le C \int_0^t \int_0^s e^{-\frac12\tilde{\nu}_{0}(t-\tau)} \left\{\int_{|u|<3N} \int_{|\xi|<5N} \left|h(\tau,X(T_1-s+\tau;T_1,X(s),u),\xi)\right|^2 \dif \xi \dif u\right\}^{\frac12} \dif \tau \dif s.
     \end{split}
    \end{equation}
    Using Lemma \ref{lem:Change of Variable}, we split the term on the RHS of \eqref{3.43} as
    \begin{equation}\label{3.44}
     \begin{split}
      &\int_0^t \int_0^s e^{-\frac12\tilde{\nu}_{0}(t-\tau)} \left\{\int_{|u|<3N} \int_{|\xi|<5N} |h(\tau,X(T_1-s+\tau;T_1,X(s),u),\xi)|^2 \dif \xi \dif u\right\}^{\frac12} \dif \tau \dif s
      \\&\le\sum_{i^1=1}^{M_1} \sum_{I^2\in\{1,\cdots,M_2\}^3} \sum_{I^3\in\{1,\cdots,M_3\}^3} \int_0^t \mathbf{1}_{\{X(s)\in\mathscr{D}_{I^2}^2\}}(s) \dif s \int_0^s e^{-\frac12\tilde{\nu}_{0}(t-\tau)} \mathbf{1}_{\mathscr{D}_{i^1}^1}(T_1-s+\tau) \dif \tau
      \\&\qquad\times \left\{\int_{|u|<3N} \int_{|\xi|<5N} |h(\tau,X(T_1-s+\tau;T_1,X(s),u),\xi)|^2 \mathbf{1}_{\mathscr{D}_{I^3}^3}(u) \dif \xi \dif u\right\}^{\frac12}
      \\&=\sum_{j=1}^3 \sum_{i^1} \sum_{I^2} \sum_{I^3} \int_0^t \mathbf{1}_{\{X(s)\in\mathscr{D}_{I^2}^2\}}(s) \dif s \int_0^s e^{-\frac12\tilde{\nu}_{0}(t-\tau)} \mathbf{1}_{\mathscr{D}_{i^1}^1}(T_1-s+\tau)
      \\&\qquad\qquad\times\mathbf{1}_{(t_{j,i^1,I^2,I^3}-\frac\varepsilon{4M^1},t_{j,i^1,I^2,I^3}+\frac\varepsilon{4M^1})}(T_1-s+\tau) \dif \tau
      \\&\qquad\times \left\{\int_{|u|<3N} \int_{|\xi|<5N} |h(\tau,X(T_1-s+\tau;T_1,X(s),u),\xi)|^2 \mathbf{1}_{\mathscr{D}_{I^3}^3}(u) \dif \xi \dif u\right\}^{\frac12}
      \\&\quad+\sum_{j=1}^3 \sum_{i^1} \sum_{I^2} \sum_{I^3} \int_0^t \mathbf{1}_{\{X(s)\in\mathscr{D}_{I^2}^2\}}(s) \dif s \int_0^s e^{-\frac12\tilde{\nu}_{0}(t-\tau)} \mathbf{1}_{\mathscr{D}_{i^1}^1}(T_1-s+\tau)
      \\&\qquad\qquad\times\left\{1-\mathbf{1}_{(t_{j,i^1,I^2,I^3}-\frac\varepsilon{4M^1},t_{j,i^1,I^2,I^3}+\frac\varepsilon{4M^1})}(T_1-s+\tau)\right\} \dif \tau
      \\&\qquad\times \left\{\int_{|u|<3N} \int_{|\xi|<5N} |h(\tau,X(T_1-s+\tau;T_1,X(s),u),\xi)|^2 \mathbf{1}_{\mathscr{D}_{I^3}^3}(u) \dif \xi \dif u\right\}^{\frac12}.
     \end{split}
    \end{equation}
    Noting
    \begin{equation}\nonumber
     \sum_{I^2\in\{1,\cdots,M_2\}^3} \int_0^t \mathbf{1}_{\{X(s)\in\mathscr{D}_{I^2}^2\}}(s) \dif s =\int_0^t \dif s,
    \end{equation}
    then the first term on the RHS of \eqref{3.44} is bounded by
    \begin{equation}\label{3.45}
     \begin{split}
      &C_N \sup_{0\le s\le t} \|h(s)\|_{L^\infty} \sum_{j=1}^3 \sum_{i^1} \sum_{I^2} \int_0^t \mathbf{1}_{\{X(s)\in\mathscr{D}_{I^2}^2\}}(s) \dif s \int_0^s e^{-\frac12\tilde{\nu}_{0}(t-\tau)} \mathbf{1}_{\mathscr{D}_{i^1}^1}(T_1-s+\tau)
      \\&\qquad\qquad\times\mathbf{1}_{(t_{j,i^1,I^2,I^3}-\frac\varepsilon{4M^1},t_{j,i^1,I^2,I^3}+\frac\varepsilon{4M^1})}(T_1-s+\tau) \dif \tau
      \\&\le C_N \sup_{0\le s\le t} \|h(s)\|_{L^\infty} \sum_{j=1}^3 \sum_{i^1} \sum_{I^2} \int_0^t \frac{\varepsilon}{M_1} e^{-\frac12\tilde{\nu}_{0}(t-s)} \mathbf{1}_{\{X(s)\in\mathscr{D}_{I^2}^2\}}(s) \dif s
      \\&\le C_N \cdot \varepsilon \sup_{0\le s\le t} \|h(s)\|_{L^\infty}.
     \end{split}
    \end{equation}
    On the other hand, it follows from \eqref{2.5} that
    \begin{equation}\nonumber
     \left|\frac{\partial X(T_1-s+\tau;T_1,X(s),u)}{\partial v}\right| \ge \delta_{*}(\varepsilon, M_1,M_2,M_3,N,T_1)>0,
    \end{equation}
    for $X(s)\in\mathscr{D}_{I^2}^2$, $u\in\mathscr{D}_{I^3}^3$ and $T_1-s+\tau \notin (t_{j,i^1,I^2,I^3}-\frac\varepsilon{4M^1},t_{j,i^1,I^2,I^3}+\frac\varepsilon{4M^1})$, which immediately yields that the second term on the RHS of \eqref{3.44} is bounded by
    \begin{equation}\label{3.46}
     \begin{split}
      &\frac{C_N}{\delta_*} \sum_{j=1}^3 \sum_{i^1} \sum_{I^2} \sum_{I^3} \int_0^t \mathbf{1}_{\{X(s)\in\mathscr{D}_{I^2}^2\}}(s) \dif s \int_0^s e^{-\frac12\tilde{\nu}_{0}(t-\tau)} \mathbf{1}_{\mathscr{D}_{i^1}^1}(T_1-s+\tau)
      \\&\qquad\qquad\times\left\{1-\mathbf{1}_{(t_{j,i^1,I^2,I^3}-\frac\varepsilon{4M^1},t_{j,i^1,I^2,I^3}+\frac\varepsilon{4M^1})}(T_1-s+\tau)\right\} \|f(\tau)\|_{L^2} \dif \tau
      \\&\le C(\varepsilon,M_1,M_2,M_3,N,T_1) \int_0^t e^{-\frac14\tilde{\nu}_{0}(t-\tau)} \|f(\tau)\|_{L^2} \dif \tau.
     \end{split}
    \end{equation}
    Substituting \eqref{3.45} and \eqref{3.46} into \eqref{3.44}, one proves
    \begin{equation}\label{3.47}
     \begin{split}
      &\int_0^t \int_0^s e^{-\frac12\tilde{\nu}_{0}(t-\tau)} \left\{\int_{|u|<3N} \int_{|\xi|<5N} |h(\tau,X(T_1-s+\tau;T_1,X(s),u),\xi)|^2 \dif \xi \dif u\right\}^{\frac12} \dif \tau \dif s
      \\&\le C_N\cdot\varepsilon \sup_{0\le s\le t} \|h(s)\|_{L^\infty}+C(\varepsilon,M_1,M_2,M_3,N,T_1) \int_0^t e^{-\frac14 \tilde{\nu}_0(t-s)} \|f(s)\|_{L^2} \dif s.
     \end{split}
    \end{equation}
    Combining \eqref{3.37}, \eqref{3.40}, \eqref{3.43} and \eqref{3.47}, one obtains that
    \begin{align}\label{est:J11}
     J_{11} &\le C \left(\frac1N+C_N\cdot\varepsilon\right) \sup_{0\le s\le t} \|h(s)\|_{L^\infty}\nonumber\\
     &\quad+C(\varepsilon,M_1,M_2,M_3,N,T_1) \int_0^t e^{-\frac14 \tilde{\nu}_0(t-s)} \|f(s)\|_{L^2} \dif s.
    \end{align}

    For $J_{12}$, it follows from Lemma \ref{lem:2.2} that
    \begin{equation}\label{ineq:J12}
     \begin{split}
      J_{12}&\le C\sup_{0\le s\le t} \|h(s)\|_{L^\infty} \int_0^t \int_0^s e^{-\frac12\tilde{\nu}_{0}(t-\tau)} \int_{\mathbb{R}^3} \left|k_{\beta}(V(s),u)\right|
      \\&\qquad\qquad\qquad\qquad\quad\times\left\{\int_{\mathbb{R}^3} (1+|\eta|)^4 \left|f(\tau,\tilde{X}(\tau),\eta)\right|^2 \dif \eta\right\}^{\frac12} \dif u \dif \tau \dif s.
     \end{split}
    \end{equation}
    We divide the estimation of \eqref{ineq:J12} into the following several cases.\\
    \noindent{\rm Case I.} For $|v|\ge N\gg1$, it follows from \eqref{est:|V(s)|} and \eqref{est:k} that
    \begin{align}\label{est:J12 |v|>=N}
      J_{12}& \le C \sup_{0\le s\le t} \|h(s)\|^2_{L^\infty} \int_0^t \int_0^s e^{-\frac12 \tilde{\nu}_0(t-\tau)} \int_{\mathbb{R}^3} \left|k_\beta(V(s),u)\right| \dif u \dif \tau \dif s
      \nonumber\\
      & \le C \sup_{0\le s\le t} \|h(s)\|^2_{L^\infty} \int_0^t \int_0^s e^{-\frac12 \tilde{\nu}_0(t-\tau)} \frac1{1+|V(s)|} \dif \tau \dif s
      \nonumber\\
      &\le \frac{C}{N} \sup_{0\le s\le t} \|h(s)\|^2_{L^\infty}.
    \end{align}

    \noindent{\rm Case II.} For either $|v|<N,|u|\ge3N$ or $|u|<3N,|\eta|\ge5N$, it follows from \eqref{est:|V(s)|} that either $|V(s)-u|\ge N,\ \forall s\in[0,t]$ or $|\eta|\ge5N$, which together with \eqref{3.38} yields that
    \begin{align}\label{est:J12 |v|<N |u|>=3N or |u|<3N |eta|>=5N}
      & C\sup_{0\le s\le t} \|h(s)\|_{L^\infty} \int_0^t \int_0^s e^{-\frac12 \tilde{\nu}_0(t-\tau)} \int_{|u|\ge3N} \left|k_{\beta}(V(s),u)\right| \left\{\int_{\mathbb{R}^3} \cdots \dif \eta\right\}^{\frac12} \dif u \dif \tau \dif s
      \nonumber\\
      &\quad+C\sup_{0\le s\le t} \|h(s)\|_{L^\infty} \int_0^t \int_0^s e^{-\frac12 \tilde{\nu}_0(t-\tau)} \int_{|u|<3N} \left|k_{\beta}(V(s),u)\right| \left\{\int_{|\eta|\ge5N} \cdots \dif \eta\right\}^{\frac12} \dif u \dif \tau \dif s
      \nonumber\\
      &\le \frac{C}{\sqrt{N}} \sup_{0\le s\le t} \|h(s)\|^2_{L^\infty}.
    \end{align}

    \noindent{\rm Case III.} $|v|<N, |u|<3N, |\eta|<5N$. This is the last remaining case. It follows from \eqref{3.41}, \eqref{3.42} and \eqref{3.47} that
    \begin{equation}\label{est:J12 |v|<N |u|<3N |eta|<5N}
     \begin{split}
      &C\sup_{0\le s\le t} \|h(s)\|_{L^\infty} \int_0^t \int_0^s e^{-\frac12 \tilde{\nu}_0(t-\tau)} \int_{|u|<3N} \left|k_{\beta}(V(s),u)\right| \left\{\int_{|\eta|<5N} \cdots \dif \eta\right\}^{\frac12} \dif u \dif \tau \dif s
      \\& \le C_N \sup_{0\le s\le t} \|h(s)\|_{L^\infty} \int_0^t \int_0^s e^{-\frac12 \tilde{\nu}_0(t-\tau)}
      \\&\qquad\qquad\qquad\times\left\{\int_{|u|<3N} \int_{|\eta|<5N} \left|h(\tau,X(T_1-s+\tau; T_1,X(s),u),\eta)\right|^2 \dif \eta \dif u \right\}^\frac12 \dif \tau \dif s
      \\&\le C(\varepsilon,M_1,M_2,M_3,N,T_1) \sup_{0\le s\le t} \|h(s)\|_{L^\infty} \int_0^t \int_0^s e^{-\frac12 \tilde{\nu}_0(t-\tau)} \|f(\tau)\|_{L^2} \dif \tau \dif s
      \\&\qquad\quad+C_N \cdot \varepsilon \sup_{0\le s\le t} \|h(s)\|^2_{L^\infty}
      \\&\le C_N \cdot \varepsilon \sup_{0\le s\le t} \|h(s)\|^2_{L^\infty}+C(\varepsilon,M_1,M_2,M_3,N,T_1) \int_0^t e^{-\frac14\tilde{\nu}_0(t-s)} \|f(s)\|^2_{L^2} \dif s.
     \end{split}
    \end{equation}
    Substituting \eqref{est:J12 |v|>=N}, \eqref{est:J12 |v|<N |u|>=3N or |u|<3N |eta|>=5N} and \eqref{est:J12 |v|<N |u|<3N |eta|<5N} into \eqref{ineq:J12}, one obtains that
    \begin{align}\label{est:J12}
     J_{12} &\le \left(\frac{C}{\sqrt{N}}+C_N \cdot \varepsilon\right) \sup_{0\le s\le t} \|h(s)\|^2_{L^\infty}\nonumber\\
     &\quad+C(\varepsilon,M_1,M_2,M_3,N,T_1) \int_0^t e^{-\frac14 \tilde{\nu}_0(t-s)} \|f(s)\|^2_{L^2} \dif s.
    \end{align}
    Substituting \eqref{est:J11} and \eqref{est:J12} into \eqref{ineq:J1}, one gets that
    \begin{equation}\label{est:J1}
     \begin{split}
      J_1 &\le CA_0^2 t e^{-\frac12\tilde{\nu}_0 t}+C A_0 \left(\frac1{\sqrt{N}}+C_N \cdot \varepsilon\right) \sup_{0\le s\le t} \left\{\|h(s)\|_{L^\infty}+ \|h(s)\|^2_{L^\infty}\right\}
      \\&\quad+ C(\varepsilon,M_1,M_2,M_3,N,T_1) \cdot A_0 \sup_{0\le s\le t} \left\{\|f(s)\|_{L^2}+\|f(s)\|^2_{L^2}\right\}.
     \end{split}
    \end{equation}

    \

    Next we estimate the nonlinear term $J_2$. Motivated by \cite{Duan-Wang}, we need to make an iteration again in the nonlinear term. In fact, it follows from \eqref{est:Gamma+}  that
    \begin{align}\label{ineq:J2}
      J_{2} &\le C \int_0^t e^{-\frac12 \tilde{\nu}_0(t-s)} \|h(s)\|_{L^\infty} \left\{\int_{\mathbb{R}^3} (1+|u|)^{-2\beta+4} \left|h(s,X(s),u)\right|^2 \dif u \right\}^\frac12 \dif s
      \nonumber\\
      &\le \frac{C}{\sqrt{N}} \sup_{0\le s\le t} \|h(s)\|_{L^\infty}^2\nonumber\\
      &\quad+C \int_0^t e^{-\frac12 \tilde{\nu}_0(t-s)} \|h(s)\|_{L^\infty} \left\{\int_{|u|\le N} (1+|u|)^{-2\beta+4} \left|h(s,X(s),u)\right|^2 \dif u \right\}^\frac12 \dif s.
    \end{align}
    Recall $[\tilde{X}(\tau),\tilde{V}(\tau)]:=[X(\tau;s,X(s),u),V(\tau;s,X(s),u)]$, we substitute the mild formulation of $h(s,X(s),u)$ into the second term on the RHS of \eqref{ineq:J2} to obtain that
    \begin{align}\label{3.56}
      &\int_0^t e^{-\frac12 \tilde{\nu}_0(t-s)} \|h(s)\|_{L^\infty} \left(\int_{|u|\le N} (1+|u|)^{-2\beta+4} \left|h(s,X(s),u)\right|^2 \dif u \right)^\frac12 \dif s
      \nonumber\\
      &\le C A_0^2 \int_0^t e^{-\frac12 \tilde{\nu}_0(t-s)} e^{-\frac12 \tilde{\nu}_0 s} \|h(s)\|_{L^\infty} \left\{\int_{|u|\le N} (1+|u|)^{-2\beta+4} \dif u \right\}^\frac12 \dif s
      \nonumber\\
      &\quad +C A_0 \int_0^t e^{-\frac12 \tilde{\nu}_0(t-s)} \|h(s)\|_{L^\infty} \Bigg\{\int_0^s e^{-\frac12\tilde{\nu}_{0}(s-\tau)} \dif \tau
      \nonumber\\
      &\qquad\qquad \times \int_{|u|\le N} (1+|u|)^{-2\beta+4}\left(\int_{\mathbb{R}^3} \left|k_{\beta}(\tilde{V}(\tau),\xi) h(\tau,\tilde{X}(\tau),\xi)\right| \dif \xi\right)^2 \dif u \Bigg\}^\frac12 \dif s
      \nonumber\\
      &\quad +C A_0 \int_0^t e^{-\frac12 \tilde{\nu}_0(t-s)} \|h(s)\|_{L^\infty} \Bigg\{\int_0^s e^{-\frac12 \tilde{\nu}_0(s-\tau)} \|h(\tau)\|_{L^\infty}^2 \dif \tau
      \nonumber\\
      &\qquad\quad \times \int_{|u|\le N} \int_{\mathbb{R}^3} (1+|u|)^{-2\beta+4} (1+|\eta|)^{-2\beta+4} \left|h(\tau,\tilde{X}(\tau),\eta)\right|^2 \dif \eta \dif u \Bigg\}^\frac12 \dif s
      \nonumber\\
      &:=J_{21}+J_{22}+J_{23}.
    \end{align}
    A direct calculation shows that
    \begin{equation}\label{est:J21}
     J_{21}\le C A_0^2 e^{-\frac12 \tilde{\nu}_0 t} \int_0^t \|h(s)\|_{L^\infty} \dif s.
    \end{equation}
    For $J_{22}$, it follows from \eqref{3.39}, \eqref{3.41} and \eqref{3.42} that
    \begin{align}\label{ineq:J22}
      J_{22}&\le C A_0 \int_0^t e^{-\frac12 \tilde{\nu}_0(t-s)} \|h(s)\|_{L^\infty} \Bigg\{\int_0^s e^{-\frac12\tilde{\nu}_{0}(s-\tau)} \dif \tau
      \nonumber\\
      &\quad \times \int_{|u|\le N} (1+|u|)^{-2\beta+4} \left(\int_{|\xi|\le3N} \left|k_{\beta}(\tilde{V}(\tau),\xi) h(s,\tilde{X}(\tau),\xi)\right| \dif \xi\right)^2 \dif u \Bigg\}^\frac12 \dif s
      \nonumber\\
      &\quad+\frac{CA_0}{N} \sup_{0\le s\le t} \|h(s)\|_{L^\infty}^2
      \nonumber\\
      &\le C A_0 \sup_{0\le s\le t} \|h(s)\|_{L^\infty} \int_0^t e^{-\frac12 \tilde{\nu}_0(t-s)}  \Bigg\{\int_0^s e^{-\frac12\tilde{\nu}_{0}(s-\tau)} \dif \tau
      \nonumber\\
      &\quad \times \int_{|u|\le N} \int_{|\xi|\le3N} (1+|u|)^{-2\beta+4} \left|h(s,X(T_1-s+\tau; T_1,X(s),u),\xi)\right|^2 \dif \xi \dif u \Bigg\}^\frac12 \dif s
      \nonumber\\
      &\quad+\frac{CA_0}{N} \sup_{0\le s\le t} \|h(s)\|_{L^\infty}^2.
    \end{align}
    By similar arguments as in \eqref{3.44}-\eqref{3.47}, one has that
    \begin{equation}\label{3.59}
     \begin{split}
      &\int_0^t e^{-\frac12 \tilde{\nu}_0(t-s)} \Bigg\{\int_0^s e^{-\frac12\tilde{\nu}_{0}(s-\tau)} \int_{|u|\le N} \int_{|\xi|\le3N} (1+|u|)^{-2\beta+4}\dif \tau
      \\&\qquad\qquad\qquad \times  \left|h(s,X(T_1-s+\tau; T_1,X(s),u),\xi)\right|^2 \dif \xi \dif u \Bigg\}^\frac12 \dif s
      \\&\le C_N \cdot \varepsilon \sup_{0\le s\le t} \|h(s)\|_{L^\infty} +C(\varepsilon,M_1,M_2,M_3,N,T_1) \left\{\int_0^t e^{-\frac12\tilde{\nu}_{0}(t-\tau)} \|f(\tau)\|^2_{L^2} \dif \tau\right\}^{\frac12},
     \end{split}
    \end{equation}
    which, together with \eqref{ineq:J22}, yields that
    \begin{align}\label{est:J22}
     J_{22} &\le A_0\left(\frac{C}{N}+C_N\cdot\varepsilon\right) \sup_{0\le s\le t} \|h(s)\|_{L^\infty}^2\nonumber\\
     &\quad+C(\varepsilon,M_1,M_2,M_3,N,T_1) \int_0^t e^{-\frac12\tilde{\nu}_{0}(t-\tau)} \|f(\tau)\|^2_{L^2} \dif \tau.
    \end{align}
    For $J_{23}$, it follows from  \eqref{3.59}  that
    \begin{align}\label{est:J23}
      J_{23}&\le CA_0 \sup_{0\le s\le t} \|h(s)\|_{L^\infty}^2 \int_0^t e^{-\frac12 \tilde{\nu}_0(t-s)} \Bigg\{\int_0^s e^{-\frac12 \tilde{\nu}_0(s-\tau)} \dif \tau \int_{|u|\le N} \int_{|\eta|\le3N}
      \nonumber\\
      &\qquad\quad \times (1+|u|)^{-2\beta+4} (1+|\eta|)^{-2\beta+4} \left|h(\tau,\tilde{X}(\tau),\eta)\right|^2 \dif \eta \dif u \Bigg\}^\frac12 \dif s
      \nonumber\\
      &\quad+\frac{CA_0}{\sqrt{N}} \sup_{0\le s\le t} \|h(s)\|_{L^\infty}^3
      \nonumber\\
      &\le A_0 \left(\frac{C}{\sqrt{N}}+C_N \cdot \varepsilon\right) \sup_{0\le s\le t} \|h(s)\|_{L^\infty}^3
      \nonumber\\
      &\qquad+A_0 C(\varepsilon,M_1,M_2,M_3,N,T_1) \left\{\int_0^t e^{-\frac12\tilde{\nu}_{0}(t-\tau)} \|f(\tau)\|^2_{L^2} \dif \tau\right\}^{\frac32}.
    \end{align}
    Combining \eqref{3.56},  \eqref{est:J21}, \eqref{est:J22}, \eqref{est:J23} with \eqref{ineq:J2}, one proves that
    \begin{equation}\nonumber
     \begin{split}
      J_2&\le C A_0^2 e^{-\frac12 \tilde{\nu}_0t} \int_0^t \|h(s)\|_{L^\infty} \dif s+A_0 \left(\frac{C}{\sqrt{N}}+C_N \cdot \varepsilon\right) \sup_{0\le s\le t} \left\{\|h(s)\|_{L^\infty}^2+\|h(s)\|_{L^\infty}^3\right\}
      \\&\quad+A_0 C(\varepsilon,M_1,M_2,M_3,N,T_1) \sup_{0\le s\le t} \left\{\|f(s)\|^2_{L^2}+\|f(s)\|^3_{L^2}\right\},
     \end{split}
    \end{equation}
    which, together with \eqref{est:J1} and \eqref{3.35}, yields that
    \begin{equation}\nonumber
     \begin{split}
      \|h(t)\|_{L^\infty}&\le \tilde{C}_3 A_0^2\left(1+\int_0^t \|h(s)\|_{L^\infty} \dif s\right) \exp\left\{-\frac14 \tilde{\nu}_0 t \right\}
      \\&\quad+A_0 \left(\frac{C}{\sqrt{N}}+C_N\cdot\varepsilon\right) \sup_{0\le s\le t} \left\{\|h(s)\|_{L^\infty}+\|h(s)\|_{L^\infty}^3\right\}
      \\&\quad+A_0 C(\varepsilon,M_1,M_2,M_3,N,T_1) \sup_{0\le s\le t} \left\{  \|f(s)\|_{L^2}+\|f(s)\|^3_{L^2}\right\},
     \end{split}
    \end{equation}
    where $\tilde{C}_3\ge1$ is some generic positive constant. Therefore the proof of Lemma \ref{lem:3.3} is completed. $\hfill\Box$
   \end{pf}

  \

  \subsection{Proof of Theorem \ref{thm:main theorem}}\label{section3.4}
   We take
   \begin{equation}\nonumber
    \tilde{C}_4=\max\left\{2,C_0,\tilde{C}_3\right\},
   \end{equation}
   then we define
   \begin{equation}\label{def:A1 and T1}
    A_1:=4\tilde{C}_4^2 A_0^2 \exp\left\{\frac{4}{\tilde{\nu}_0} \tilde{C}_4 A_0^2\right\} \quad\text{and}\quad T_1:=\frac{8}{\tilde{\nu}_0} \left(\ln{A_1}+|\ln\delta|\right),
   \end{equation}
   where $C_0>$ and $\delta>0$ are the constants introduced in Theorem \ref{thm:kim}. We emphasize that the above $A_1$ in \eqref{def:A1 and T1} depends only on $A_0$, and $T_1$ depends only on $\delta$ and $A_0$.

   Assume that $\|f_0\|_{L^2}\le\kappa_1(A_1,T_1)$ where $\kappa_1$ is defined in Lemma \ref{lem:3.3}. Hence, it follows from \eqref{a priori assumption}, \eqref{def:A1 and T1} and Lemmas \ref{lem:3.2} and \ref{lem:3.3} that
   \begin{equation}\label{3.63}
    \|h(t)\|_{L^\infty}\le \tilde{C}_4 A_0^2 \left(1+\int_0^t\|h(s)\|_{L^\infty} \dif s \right) e^{-\f14\tilde{\nu}_0t}+D,\quad t\in[0,T_1],
   \end{equation}
   where
   \begin{equation}\label{3.64}
    \begin{split}
     D&:=A_0\left(\frac{\tilde{C}_3}{\sqrt{N}}+C_N\cdot\varepsilon\right) \left\{A_1+A_1^3\right\}
     \\
     &\quad+A_0 C(\varepsilon,M_1,M_2,M_3,N,T_1) \tilde{C}_1^3 \left\{\|f_0\|_{L^2} e^{\tilde{C}_1 A_1 T_1}+\left(\|f_0\|_{L^2} e^{\tilde{C}_1 A_1 T_1}\right)^3\right\}.
    \end{split}
   \end{equation}

   Now we define
   \begin{equation}\nonumber
    H(t):=1+\int_0^t \|h(s)\|_{L^\infty} \dif s,
   \end{equation}
   then, \eqref{3.63} is rewritten as
   \begin{equation}\nonumber
    H'(t)\le \tilde{C}_4 A_0^2 e^{-\f14\tilde{\nu}_0t} H(t)+D.
   \end{equation}
   Hence it holds that
   \begin{equation}\nonumber
    \frac{\dif}{\dif t} \left(H(t)\exp{\left\{-\frac{4}{\tilde{\nu}_0} \tilde{C}_4 A_0^2\left(1-e^{-\frac14\tilde{\nu}_0 t}\right) \right\}}\right) \le D,
   \end{equation}
   which yields immediately that
   \begin{equation}\label{3.65}
    H(t)\le (1+Dt)\exp\left\{\frac{4}{\tilde{\nu}_0} \tilde{C}_4 A_0^2\right\},\quad\forall t\in[0,T_1].
   \end{equation}
   Substituting \eqref{3.65} back into \eqref{3.63}, one has
   \begin{equation}\label{3.66}
    \begin{split}
     \|h(t)\|_{L^\infty}&\le \tilde{C}_4 A_0^2 \exp\left\{\frac{4}{\tilde{\nu}_0} \tilde{C}_4 A_0^2\right\} (1+Dt) e^{-\f14\tilde{\nu}_0t}+D
     \\&\le \frac{1}{4\tilde{C}_4} A_1\left(1+\frac{8}{\tilde{\nu}_0}D\right)e^{-\f18\tilde{\nu}_0t}+D.
    \end{split}
   \end{equation}
   Noting \eqref{def:A1 and T1} and \eqref{3.64}, we firstly choose $N>0$ large enough, then $\varepsilon>0$ sufficiently small, and finally let $\|f_0\|_{L^2}\le\kappa_2$ with $\kappa_2=\kappa_2(\delta, A_0)>0$ further sufficiently small, such that
   \begin{equation}\nonumber
    D\le \min\left\{\f{\tilde{\nu}_0}{32}, \frac{\delta}{8}\right\},
   \end{equation}
   which, together with \eqref{3.66} yields that
   \begin{equation}\label{est:uniform estimate}
    \|h(t)\|_{L^\infty} \le\frac{5}{16\tilde{C}_4} A_1e^{-\frac18 \tilde{\nu}_0t}+D\le \frac{1}{2\tilde{C}_4} A_1,
   \end{equation}
   for all $t\in[0,T_1]$. Hence, we have closed the {\it a priori} assumption \eqref{a priori assumption} over $t\in[0,T_1]$ provided that
   \begin{equation}\nonumber
    \|f_0\|_{L^2}\le\kappa_0:=\min\{\kappa_1,\kappa_2\}.
   \end{equation}
   Note that $\kappa_0>0$ depends only on $\delta$ and $A_0$.

   Using the uniform estimate \eqref{est:uniform estimate} and the local existence Theorem \ref{thm:local existence}, we can extend the Boltzmann solution to the time interval $t\in[0,T_1]$, see \cite{Duan-Wang} for more details. Next, for the case  $t\ge T_1$, we note from the first inequality of \eqref{est:uniform estimate} and the definition \eqref{def:A1 and T1} for $T_1$ that
   \begin{equation}\nonumber
    \|h(T_1)\|_{L^\infty} \le \frac5{16\tilde{C}_4} A_1 \exp\left\{-\frac{\tilde{\nu}_0}{8}T_1\right\}+\frac{\delta}{8}\le \frac{5\delta}{16\tilde{C}_4}+\frac{\delta}{8}<\frac{1}{2}\delta.
   \end{equation}
   With $t=T_1$ as the initial time and applying Theorem \ref{thm:kim}, we can extend the Boltzmann solution $f(t)$ from $[0,T_1]$ to $[0,\infty)$, and thus obtain the unique solution $f(t)$ globally in time on $[0,\infty)$ such that $F(t,x,v)=\mu_{E}+\sqrt{\mu_E}f(t,x,v)\ge0$ and $\sup_{t\ge0} \|w_\beta f(t)\|_{L^\infty}\le A_1$. This proves the global existence and uniqueness of solutions in weighted $L^\infty$ space.

   For the large time behavior of the obtained solution, we note that as an immediate consequence of Theorem \ref{thm:kim}, it holds that
   \begin{equation}\label{3.68}
    \|h(t)\|_{L^\infty}\le C_0 \|h(T_1)\|_{L^\infty} e^{-\lambda_0 (t-T_1)}\le C_0 \delta e^{-\lambda_0 (t-T_1)},
   \end{equation}
   for all $t\ge T_1$. By taking
   \begin{equation}\nonumber
    \tilde{C}_0:=4\tilde{C}_4^3 \quad\text{and}\quad \tilde{\lambda}_0:=\min\left\{\lambda_0,\f{\tilde{\nu}_0}{8}\right\},
   \end{equation}
   it follows from \eqref{est:uniform estimate}, \eqref{3.68} and a direct computation that
   \begin{equation}\nonumber
    \begin{split}
     \|h(t)\|_{L^\infty}&\le \max\left\{\frac12,C_0\right\} A_1 e^{-\tilde{\lambda} t}\le \tilde{C}_4 A_1 e^{-\tilde\lambda_0 t}
    \\&=4\tilde{C}_4^3 A_0^2 \exp\left\{\frac{4\tilde{C}_4}{\tilde{\nu}_0}A_0^2\right\} e^{-\tilde\lambda_0 t}
    \\&= \tilde{C}_0 A_0^2\exp\left\{\frac{\tilde{C}_0}{\tilde{\nu}_0}A_0^2\right\} e^{-\tilde\lambda_0 t}.
    \end{split}
   \end{equation}
   Therefore the proof of Theorem \ref{thm:main theorem} is completed. $\hfill\Box$

 \

 \section{Appendix}
  The following theorem is devoted to the existence of local solution to Boltzmann equation with large external potential and large $L^\infty$ initial data.
  \begin{theorem}[Local Existence]\label{thm:local existence}
   Let $0\le\gamma\le1$. Assume \eqref{def:the C3 norm of Phi}, and
   \begin{equation}\nonumber
    F_0(x,v)=\mu_E(x,v)+\sqrt{\mu_E(x,v)}f_0(x,v) \ge 0,
   \end{equation}
   with ${\left\|w_\beta f_0 \right\|}_{L^\infty}<+\infty$, then there exists a positive time
   \begin{equation}\nonumber
    T_0:=\left(8C \left[1+{\left\|w_\beta f_0\right\|}_{L^\infty}\right]\right)^{-1}>0,
   \end{equation}
   such that the Boltzmann equation \eqref{eq:main} admits a unique solution
   \begin{equation}\nonumber
    F(t,x,v)=\mu_E(x,v)+\sqrt{\mu_E(x,v)}f(t,x,v)\ge0
   \end{equation}
   satisfying
   \begin{equation}\nonumber
    {\left\|w_\beta f(t) \right\|}_{L^\infty} \le 2{\left\|w_\beta f_0 \right\|}_{L^\infty},\ \text{for}\ 0\le t \le T_0,
   \end{equation}
   where $C\ge1$ is some positive constant depending only on $\beta,\alpha,\gamma,M$. In addition, the conservations of mass \eqref{eq:conservation of mass}, energy \eqref{eq:conservation of energy}, and degenerate momentum \eqref{eq:degenerate conservation of momentum} as well as the additional entropy inequality \eqref{ineq:entropy inequality} hold.
  \end{theorem}
  \begin{remark}
   In fact, the above local existence result can be obtained by making a slight modification to the proof of Proposition 2.1 in \cite{Duan-Huang-Wang-Yang}. For completeness, we put its proof in the following.
  \end{remark}

  \noindent{\bf Proof of Theorem \ref{thm:local existence}.}
To prove the local existence for the Boltzmann equation with large external potential and large initial data, we consider the following iteration, for $n=0,1,2,\dots$,
   \begin{equation}\label{eq:iteration}
    \begin{cases}
     \partial_t F^{n+1}+v \cdot \nabla_xF^{n+1} -\nabla\Phi(x) \cdot \nabla_vF^{n+1} +Q_-(F^n,F^{n+1})=Q_+(F^n,F^n),
     \\[2mm] F^{n+1}(0,x,v)=F_0(x,v) \ge 0,
    \end{cases}
   \end{equation}
   where $F^0(t,x,v)=\mu_E(x,v)$. Along the characteristic line \eqref{eq:characteristics}, the mild form of \eqref{eq:iteration} is
   \begin{equation}\nonumber
    \begin{split}
     F^{n+1}(t,x,v) &=e^{-\int_0^t g_1^n(\tau) \dif \tau}F_0(X(0),V(0))
     \\&\quad+\int_0^t e^{-\int_s^t g_1^n(\tau) \dif \tau} Q_+(F_n,F_n)(s,X(s),V(s)) \dif s,
    \end{split}
   \end{equation}
   where
   \begin{equation}\nonumber
    g_1^n(\tau)=R(F^n)(\tau,X(\tau),V(\tau)),
   \end{equation}
   with $R(F)$ defined in \eqref{def:R(F)}. By the induction argument, it is easy to show that
   \begin{equation}\label{est:Fn>=0}
    F^n(t,x,v)\ge0, \quad \forall n\ge0.
   \end{equation}
   To take the limit $n\rightarrow\infty$, one needs to obtain some uniform estimate. In fact, to show the uniform estimate for the approximation sequence $F^{n}$, we turn to estimate $h^n:=w_\beta f^n$ with
   \begin{equation}\nonumber
    f^{n}:=\frac{F^{n}-\mu_E}{\sqrt{\mu_E}} \quad \text{and} \quad f^n|_{t=0}=f_0:= \frac{F_0-\mu_E}{\sqrt{\mu_E}},\quad \forall n\ge0,
   \end{equation}
   then the equation \eqref{eq:iteration} is rewritten as
   \begin{equation}\nonumber
    \begin{split}
     &\frac{\partial h^{n+1}}{\partial t} +v \cdot \nabla_x h^{n+1} -\nabla\Phi \cdot \nabla_v h^{n+1} + R(F^n)h^{n+1}
     \\&=e^{-\frac\Phi2} w_\beta \Gamma_+\left(\frac{\sqrt\mu h^n}{w_\beta},\frac{\sqrt\mu h^n}{w_\beta}\right)+e^{-\Phi} w_\beta K\left(\frac{h^n}{w_\beta}\right).
    \end{split}
   \end{equation}
   Integrate along the characteristics \eqref{eq:characteristics} to obtain
   \begin{equation}\label{eq:hn+1}
    \begin{split}
     &h^{n+1}(t,x,v)
     \\&=e^{-\int_0^t g_1^n(\tau) \dif \tau} h_0(X(0),V(0))
     \\&\quad+\int_0^t e^{-\int_s^t g_1^n(\tau) \dif \tau-\Phi(X(s))} w_\beta(X(s),V(s)) K\left(\frac{h^n}{w_\beta}\right)(s,X(s),V(s)) \dif s
     \\& \quad+\int_0^t e^{-\int_s^t g_1^n(\tau) \dif \tau-\frac{\Phi(X(s))}2} \frac{w_\beta(X(s),V(s))}{\sqrt{\mu(V(s))}} Q_+\left(\frac{\sqrt\mu h^n}{w_\beta},\frac{\sqrt\mu h^n}{w_\beta}\right)(s,X(s),V(s)) \dif s
     \\&:=H_1+H_2+H_3.
    \end{split}
   \end{equation}
   Since $g_1(\tau)=R(F^n)(\tau,X(\tau),V(\tau))\ge0$, one has that
   \begin{equation}\label{est:H1}
    |H_1| \le \|h_0\|_{L^\infty}.
   \end{equation}
   For $H_2$, it follows from \eqref{est:k} that
   \begin{equation}\label{est:H2}
    |H_2|\le C \int_0^t \|h^{n}(s)\|_{L^\infty} \dif s.
   \end{equation}
   For $H_3$, it follows from \eqref{est:Gamma+ alpha>=2} that
   \begin{equation}\nonumber
    \frac{w_\beta(X(s),V(s))}{\sqrt{\mu(V(s))}} \left| Q_+\left(\frac{\sqrt\mu h^n}{w_\beta},\frac{\sqrt\mu h^n}{w_\beta}\right)(s,X(s),V(s))\right| \le C\|h^n(s)\|_{L^\infty}^2,
   \end{equation}
   which immediately yields that
   \begin{equation}\label{est:H3}
    |H_3| \le C\int_0^t \|h^n(s)\|_{L^\infty}^2 \dif s.
   \end{equation}
   Substituting \eqref{est:H1}, \eqref{est:H2} and \eqref{est:H3} into \eqref{eq:hn+1}, one obtains that
   \begin{equation}\nonumber
    \|h^{n+1}(t)\|_{L^\infty} \le \|h_0\|_{L^\infty}+Ct\left\{\sup_{0\le s\le t} \|h^n(s)\|_{L^\infty}+\sup_{0\le s\le t}\|h^n(s)\|_{L^\infty}^2 \right\},
   \end{equation}
   where $C\ge1$ depends only on $\beta,\gamma$ and $M$. Take
   \begin{equation}\nonumber
    T_0=\left(4C\left[1+\|h_0\|_{L^\infty}\right]\right)^{-1}<1,
   \end{equation}
   and by induction arguments, we get the uniform estimate
   \begin{equation}\label{est:uniform estimate for hn}
    \|h^n(t)\|_{L^\infty} \le 2\|h_0\|_{L^\infty}, \quad \forall t\in[0,T_0],
   \end{equation}
   which yields that $\left\{\frac{h^n}{\sqrt{w_\beta}}\right\}_{n=1}^\infty$ is a sequence in $L^\infty([0,T_0]\times\mathbb{T}^3\times\mathbb{R}^3)$.

   It is noted that the approximation sequence $h^n$ itself does not converge in $L^\infty([0,T_0]\times\mathbb{T}^3\times\mathbb{R}^3)$. However, we can prove that $\left\{ \frac{h^n}{\sqrt{w_\beta}}\right\}_{n=1}^\infty$ is a Cauchy sequence in $L^\infty([0,T_0]\times\mathbb{T}^3\times\mathbb{R}^3)$. In fact, it follows from \eqref{eq:hn+1} that
   \begin{equation}\label{ineq:hn+1-hn}
    \begin{split}
     &\left|(h^{n+1}-h^n)(t,x,v)\right|
     \\& \le |h_0(X(0),V(0))| \int_0^t \left|(g_2^n-g_2^{n-1})(\tau)\right| \dif \tau
     \\& \quad+\int_0^t w_\beta(X(s),V(s)) \left|K\left(\frac{h^n}{\tilde{w}}\right)(s,X(s),V(s))\right| \int_s^t \left|(g_2^n-g_2^{n-1})(\tau)\right| \dif \tau \dif s
     \\& \quad+\int_0^t w_\beta(X(s),V(s)) \left|K\left(\frac{h^n-h^{n-1}}{w_\beta}\right)(s,X(s),V(s))\right| \dif s
     \\& \quad+\int_0^t \frac{w_\beta(X(s),V(s))}{\sqrt{\mu(V(s))}} \left|Q_+\left(\frac{\sqrt\mu h^n}{w_\beta},\frac{\sqrt\mu h^n}{w_\beta}\right)(s,X(s),V(s))\right| \int_s^t \left|(g_2^n-g_2^{n-1})(\tau)\right| \dif \tau \dif s
     \\& \quad+\int_0^t \frac{w_\beta(X(s),V(s))}{\sqrt{\mu(V(s))}} \left|Q_+\left(\frac{\sqrt\mu h^n}{w_\beta},\frac{\sqrt\mu (h^n-h^{n-1})}{w_\beta}\right)(s,X(s),V(s))\right| \dif s
     \\& \quad+\int_0^t \frac{w_\beta(X(s),V(s))}{\sqrt{\mu(V(s))}} \left|Q_+\left(\frac{\sqrt\mu (h^n-h^{n-1})}{w_\beta},\frac{\sqrt\mu h^{n-1}}{w_\beta}\right)(s,X(s),V(s)) \right|\dif s
     \\&:=H_4+H_5+H_6+H_7+H_8+H_9.
    \end{split}
   \end{equation}
   By a direct calculation, we have
   \begin{equation}\label{est:g2n-g2n-1}
    \left|(g_2^{n}-g_2^{n-1})(\tau)\right| \le C \nu(V(\tau)) \left\| \left(\frac{h^n}{\sqrt{w_\beta}}-\frac{h^{n-1}}{\sqrt{w_\beta}} \right)(\tau)\right\|_{L^\infty}.
   \end{equation}
   It follows from \eqref{est:|V(s)|} that
   \begin{equation}\label{4.10}
    \int_s^t |\nu(V(\tau))| \dif \tau \le Ct(1+|v|^2)^\frac\gamma2 \le Ct \sqrt{w_\beta(x,v)}
   \end{equation}
   Using \eqref{est:g2n-g2n-1} and \eqref{4.10}, a direct calculation shows that
   \begin{equation}\label{est:H4+H5+H7}
    H_4+H_5+H_7 \le Ct \sqrt{w_\beta(x,v)} \|h_0\|_{L^\infty} \cdot\sup_{0\le s\le t} \left\| \left(\frac{h^n}{\sqrt{w_\beta}}-\frac{h^{n-1}}{\sqrt{w_\beta}} \right)(s)\right\|_{L^\infty}.
   \end{equation}

   Again, \eqref{est:|V(s)|} implies that
   \begin{equation}\label{4.12}
    w_\beta(X(s),V(s))\le C w_\beta(x,v).
   \end{equation}
   Therefore, it follows from \eqref{est:Gamma+ alpha>=2} that
   \begin{equation}\label{est:H6}
    \begin{split}
     |H_6| &\le C \int_0^t w_\beta(X(s),V(s)) \left|K\left(\frac{h^n-h^{n-1}}{w_\beta}\right)(s,X(s),V(s))\right| \dif s
     \\& \le C \int_0^t \sqrt{w_\beta(X(s),V(s))} \left\| \left(\frac{h^n}{\sqrt{w_\beta}}-\frac{h^{n-1}}{\sqrt{w_\beta}}\right)(s) \right\|_{L^\infty} \int_{\mathbb{R}^3} |k(V(s),u)| \sqrt{\frac{w_\beta(X(s),V(s))}{w_\beta(X(s),u)}} \dif u \dif s
     \\& \le C t \sqrt{w_\beta(x,v)} \sup_{0\le s\le t} \left\| \left(\frac{h^n}{\sqrt{w_\beta}}-\frac{h^{n-1}}{\sqrt{w_\beta}}\right)(s) \right\|_{L^\infty},
    \end{split}
   \end{equation}

   For $H_8$, we note from \eqref{4.12} that
   \begin{equation}\nonumber
    \begin{split}
     &\frac{w_\beta(X(s),V(s))}{\sqrt{\mu(V(s))}} \left|Q_+\left(\frac{\sqrt\mu h^n}{w_\beta},\frac{\sqrt\mu (h^n-h^{n-1})}{w_\beta}\right)(s,X(s),V(s))\right|
     \\& \le C \sqrt{w_\beta(x,v)} \|h^n(s)\|_{L^\infty} \left\| \left(\frac{h^n}{\sqrt{w_\beta}}-\frac{h^{n-1}}{\sqrt{w_\beta}} \right)(s)\right\|_{L^\infty}
     \\& \quad \times \int_{\mathbb{R}^3} \int_{\mathbb{S}^2} B(V(s)-u,\omega) \sqrt{\mu(u)} \frac1{w_\beta(X(s),u')} \sqrt{\frac{w_\beta(X(s),V(s))}{w_\beta(X(s),v')}} \dif \omega \dif u
     \\&\le C\nu(v)\|h^n(s)\|_{L^\infty} \left\| \left(\frac{h^n}{\sqrt{w_\beta}}-\frac{h^{n-1}}{\sqrt{w_\beta}} \right)(s)\right\|_{L^\infty}
     \\& \le C\sqrt{w_\beta(x,v)} \|h^n(s)\|_{L^\infty} \left\| \left(\frac{h^n}{\sqrt{w_\beta}}-\frac{h^{n-1}}{\sqrt{w_\beta}} \right)(s)\right\|_{L^\infty},
    \end{split}
   \end{equation}
   which yields that
   \begin{equation}\label{est:H8}
    H_8 \le Ct \sqrt{w_\beta(x,v)} \|h_0\|_{L^\infty} \sup_{0\le s\le t} \left\| \left(\frac{h^n}{\sqrt{w_\beta}}-\frac{h^{n-1}}{\sqrt{w_\beta}} \right)(s)\right\|_{L^\infty}.
   \end{equation}
   Similarly, we have
   \begin{equation}\label{est:H9}
    H_9 \le Ct \sqrt{w_\beta(x,v)} \|h_0\|_{L^\infty} \sup_{0\le s\le t} \left\| \left(\frac{h^n}{\sqrt{w_\beta}}-\frac{h^{n-1}}{\sqrt{w_\beta}} \right)(s)\right\|_{L^\infty}.
   \end{equation}
   Substituting \eqref{est:H4+H5+H7},\eqref{est:H6}, \eqref{est:H8},\eqref{est:H9} into \eqref{ineq:hn+1-hn}, one obtains that
   \begin{equation}\nonumber
    \begin{split}
     \sup_{0\le t\le T_0} \left\|\left(\frac{h^{n+1}}{\sqrt{w_\beta}}-\frac{h^n}{\sqrt{w_\beta}} \right)(t)\right\|_{L^\infty} & \le C T_0 \left(\|h_0\|_{L^\infty}+1\right) \sup_{0\le t\le T_0} \left\| \left(\frac{h^n}{\sqrt{w_\beta}}-\frac{h^{n-1}}{\sqrt{w_\beta}} \right)(t)\right\|_{L^\infty}
     \\& \le \frac12 \sup_{0\le t\le T_0} \left\|\left(\frac{h^n}{\sqrt{w_\beta}}-\frac{h^{n-1}}{\sqrt{w_\beta}} \right)(t)\right\|_{L^\infty}.
    \end{split}
   \end{equation}
   This inequality shows that $\sqrt{w_\beta}f^n$ or equivalently $\frac{h^n}{\sqrt{w_\beta}}$ is a Cauchy sequence in $L^\infty([0,T_0]\times\mathbb{T}^3\times\mathbb{R}^3)$. Therefore, there exists a limit function $f$ such that
   \begin{equation}\nonumber
    \sup_{0\le t\le T_0} \left\|\left(\sqrt{w_\beta}f^n-\sqrt{w_\beta}f \right)(t)\right\|_{L^\infty} \rightarrow 0,\ \text{as}\ n\rightarrow\infty,
   \end{equation}
   and the limit function $F:=\mu_E(x,v)+\sqrt{\mu_E(x,v)}f(t,x,v)$ is indeed a mild solution to the Boltzmann equation \eqref{eq:main}. Moreover, it follows from \eqref{est:Fn>=0} and \eqref{est:uniform estimate for hn} that
   \begin{equation}\nonumber
    \begin{cases}
     F(t,x,v)=\mu_E(x,v)+\sqrt{\mu_E(x,v)}f(t,x,v)\ge0,
     \\[2mm] \sup_{0\le t\le T_0} \left\|w_\beta f(t)\right\|_{L^\infty} \le 2\left\|w_\beta f_0\right\|_{L^\infty}.
    \end{cases}
   \end{equation}

   For the uniqueness. Let $\tilde{F}:=\mu_E(x,v)+\sqrt{\mu_E(x,v)}\tilde{f}(t,x,v)$ be another mild solution to the Boltzmann equation \eqref{eq:main} with initial data \eqref{initial data} and satisfy
   \begin{equation}\nonumber
    \sup_{0\le t\le T_0} \|\tilde{h}(t)\|_{L^\infty}<+\infty,
   \end{equation}
   with $\tilde{h}:=w_\beta(x,v)\tilde{f}$. By the same argument as in \eqref{ineq:hn+1-hn}, we deduce
   \begin{equation}\nonumber
    \left\|\left(\frac{h}{\sqrt{w_\beta}}-\frac{\tilde{h}}{\sqrt{w_\beta}}\right)(t)\right\|_{L^\infty} \le C\left(\|h_0\|_{L^\infty}+1\right) \int_0^t \left\|\left(\frac{h}{\sqrt{w_\beta}}-\frac{\tilde{h}}{\sqrt{w_\beta}}\right)(s)\right\|_{L^\infty} \dif s
   \end{equation}
   Hence uniqueness follows immediately from the Gronwall's inequality.

   Finally, multiplying both sides of \eqref{eq:iteration} by $1,v,\frac{|v|^2}2+\Phi(x)$ and $F^{n+1}$, integrating by parts and taking the limit $n\rightarrow\infty$, we obtain the corresponding result of the conservation of mass, energy and degenerate momentum. $\hfill\Box$

   \

   \noindent{\bf Acknowledgments.} Yong Wang is partially supported by National Natural Sciences Foundation of China No. 11771429 and 11688101.


\begin{thebibliography}{99}
  \bibitem{Asano1}
  \newblock K. Asano,
  \emph{Almost transversality theorem in the classical dynamical system},
  \newblock J. Math. Kyoto Univ. 34 (1994), 87-94.

  \bibitem{Asano2}
  \newblock K. Asano,
  \emph{Local solutions to the initial boundary value problem for the Boltzmann equation with an external force},
  \newblock J. Math. Kyoto Univ. 24 (1984), 225-238.

  \bibitem{BGGL}
 \newblock  C. Bardos, I.~M. Gamba, F. Golse and C. Levermore,
 \emph{Global solutions of the Boltzmann equation over near global Maxwellians with small mass,}
 \newblock Comm. Math. Phys., 346 (2016), no. 2, 435-467.

 \bibitem{D-Lion}
 R.J. DiPerna, P.-L. Lions, On the Cauchy problem for Boltzmann equation: Global existence and weak stability,
 {\it Ann. of Math.} \textbf{130} (1989), 321--366.
	
  \bibitem{Drange}
  \newblock H. Drange,
  \emph{On the Boltzmann equation with external forces},
  \newblock SIAM J. Appl. Math. 34 (1978), 577-592.

  \bibitem{Duan}
  \newblock R.J. Duan,
  \emph{Stability of the Boltzamnn equation with potential forces on torus},
  \newblock Phys. D,  238 (2009), 1808-1820

  \bibitem{Duan-Huang-Wang-Yang}
  \newblock R.J. Duan, F.M. Huang, Y. Wang, T. Yang,
  \emph{Global well-posedness of the Boltzmann equation with large amplitude initial data},
  \newblock Arch. Rational. Mech. Anal. 225 (2017), 375-424.

\bibitem{DHWZ}
 \newblock R.J. Duan, F.M. Huang, Y. Wang, Z. Zhang,
\emph{Effects of soft interaction and non-isothermal boundary upon long-time dynamics of rarefied gas},
 \newblock  arXiv:1807.05700
	
  \bibitem{Duan-Strain}
  \newblock R.J. Duan, R. Strain,
  \emph{Optimal time decay of the Vlasov-Poisson-Boltzmann system in $\Bbb R^3$},
  \newblock Arch. Rational. Mech. Anal. 199 (2011), 291-328.

  \bibitem{Duan-Ukai-Yang-Zhao}
  \newblock R.J. Duan, S. Ukai, T. Yang, H. Zhao,
  \emph{Optimal decay estimates on the linearized Boltzmann equation with time dependent force and their applications},
  \newblock Comm. Math. Phys., 277 (2008), no. 1, 189-236.

    \bibitem{Duan-Wang}
  \newblock R.J. Duan, Y. Wang
  \emph{The Boltzmann Equation with Large-amplitude Initial Data in Bounded Domains},
  \newblock  arXiv:1703.07978

\bibitem{DuanYZ}
  \newblock R.J. Duan,  T. Yang, C.J. Zhu,
   \emph{Boltzmann equation with external force and Vlasov-Poisson-Boltzmann system in infinite vacuum.}
\newblock Discrete Contin. Dyn. Syst. 16 (2006), no. 1, 253-277.

\bibitem{DuanYZ-1}
\newblock R.J. Duan,  T. Yang, C.J. Zhu,
\emph{Global existence to Boltzmann equation with external force in infinite vacuum.}
\newblock J. Math. Phys. 46 (2005), no. 5, 053307, 13 pp.

  \bibitem{Glassey}
  \newblock R.T. Glassey,
  \emph{The cauchy problem in kinetic theory,}
  \newblock Society for Industrial and Applied Mathematics (SIAM), Philadelphia, 1996.

  \bibitem{Guo-04}  Y. Guo,  The Vlasov-Poisson-Boltzmann system near vacuum. Comm. Math. Phys. 218 (2001), 293–313.

  \bibitem{Guo-03}
  Y. Guo, 
  {Classical solutions to the Boltzmann equation for molecules with an angular cutoff},
  {\it Arch. Ration. Mech. Anal.} {\bf 169} (2003), no. 4, 305--353.


  \bibitem{Guo1}
  \newblock Y. Guo,
  \emph{Decay and continuity of the Boltzmann equation in Bounded domains},
  \newblock  Arch. Rational. Mech. Anal. 197 (2010), 713-809.

  \bibitem{Guo2}
  \newblock Y. Guo,
  \emph{Bounded solutions for the Boltzmann equation},
  \newblock  Quart. Appl. Math. 68 (2010), no.1, 143-148.

  \bibitem{IS}
  R. Illner and M. Shinbrot,
  {Global existence for a rare gas in an infinite vacuum}, {\it Comm. Math. Phys.} \textbf{95} (1984), 117--126.


  \bibitem{KS}
  S. Kaniel and M. Shinbrot,
  The Boltzmann equation. I. Uniqueness and local existence,
  {\it Comm. Math. Phys.} {\bf 58} (1978), no. 1, 65-84.

  \bibitem{Kim}
  \newblock C. Kim,
  \emph{Boltzmann equation with a large potential in a periodic box},
  \newblock Comm. Partial Diff. Eqns.  39 (2014), 1393-1423.

  \bibitem{Li-Yu}
  \newblock F.C. Li,  H.J. Yu,
  \emph{Global existence of classical solutions to the Boltzmann equationwith external force for hard potentials},
  \newblock Int. Math. Res. Not. 2008. https://doi.org/10.1093/imrn/rnn112

  \bibitem{LYY}
  \newblock T. Liu, T. Yang, and  S. H. Yu,
 \emph{Energy method for the Boltzmann equation,}
  {\it Physica D} \textbf{188} (2004), 178--192.


\bibitem{Tabata}
\newblock M. Tabata,
\emph{Decay of solutions to the mixed problem with the periodicity boundary condition for the linearized Boltzmann equation with conservative external force.}
\newblock Comm. Part. Diff. Eqs. 18 (1993), 1823-1846.

  \bibitem{Uk}
  S. Ukai,
  {On the existence of global solutions of mixed problem for non-linear Boltzmann equation}, {\it Proc. Jpn. Acad.} {\bf 50} (1974), 179--184.



  \bibitem{UA}
  {S. Ukai and K. Asano, On the Cauchy problem of the Boltzmann equation with a soft potential, {\it Publ. Res. Inst. Math. Sci.} {\bf 18} (1982), no. 2, 477--519.}


  \bibitem{Ukai-Yang-Zhao}
  \newblock S. Ukai, T. Yang, H.J.  Zhao,
  \emph{Global solution to the Boltzmann equation with external forces},
  \newblock  Anal. Appl. (Singap.) 3 (2005), no. 2, 157-193.

\bibitem{UY-2006}  S. Ukai, T. Yang,
\emph{The Boltzmann equation in the space $L^2\cap L^\infty_{x,v}\beta$: Global and time-periodic solutions. }
\newblock Anal. Appl. (Singap.) 4 (2006), no. 3, 263-310.

\bibitem{Yu}
\newblock H.J. Yu,
\emph{Global classsical solutions to the Boltzmann equation with external force},
\newblock Commun. Pure Appl. Anal. 8 (2009), 1647-1668.
 \end{thebibliography}
\end{document}